%% file: ms.tex
\documentclass[3p]{elsarticle}

\let\today\relax 
\makeatletter
\def\ps@pprintTitle{%
    \let\@oddhead\@empty
    \let\@evenhead\@empty
    \def\@oddfoot{\footnotesize\itshape
         {Submitted preprint} \hfill\today}%
    \let\@evenfoot\@oddfoot
    }
\makeatother

\usepackage[utf8]{inputenc}

\usepackage{amsfonts}
\usepackage{amsmath}
\usepackage{amssymb}
\usepackage{empheq}
\usepackage{mathrsfs} 
\usepackage{dsfont} 
\usepackage{graphicx}
\usepackage{enumitem}

\usepackage{graphicx}
\usepackage{epstopdf}
\usepackage{algorithmic} 
\usepackage{subcaption}
\usepackage{float}

\usepackage{bm}
\usepackage{tikz}
\usepackage{tkz-euclide}
\usepackage{csvsimple}

\usepackage{pgfplots}
\usepackage{floatflt}
\usepackage{booktabs}

\usepackage{wrapfig}
\usepackage{cleveref}

\pgfplotsset{compat=1.13}

\usetikzlibrary{plotmarks}
\usetikzlibrary{arrows}
\graphicspath{{./Graphics/}}

\input{ALL_ABBR.tex}

\ifpdf
  \DeclareGraphicsExtensions{.eps,.pdf,.png,.jpg}
\else
  \DeclareGraphicsExtensions{.eps}
\fi

\title{
Gradient-based shape optimization for the reduction of particle erosion in bended pipes
}

\author[1]{Raphael Hohmann\corref{cor1}%
}
\ead{raphael.hohmann@itwm.fraunhofer.de}

\author[1]{Christian Leith\"auser 
}
\ead{christian.leithaeuser@itwm.fraunhofer.de}

\cortext[cor1]{Corresponding author}
\address[1]{Fraunhofer ITWM, Fraunhofer-Platz 1, 67663 Kaiserslautern, Germany}

\newtheorem{thm}{Theorem}
\newtheorem{lem}[thm]{Lemma}
\newtheorem{defn}{Definition}
\newdefinition{rmk}{Remark}
\newproof{pf}{Proof}
\newtheorem{alg}{Algorithm}

\begin{document}

\begin{abstract}
  \input{Abstract.tex}
\end{abstract}

\begin{keyword}
  shape optimization \sep particle erosion \sep partial differential equations \sep Dean vortices
\end{keyword}

\maketitle

\pgfplotsset{compat=newest}
\pgfplotsset{plot coordinates/math parser=false}
\newlength\figureheight
\newlength\figurewidth

%
%
%
%
\section{Introduction}
\label{sec:Introduction}
\input{Introduction.tex}

\section{Problem formulation}
\label{sec:ProblemFormulation}
\input{ProblemFormulation.tex}

\section{Shape calculus}
\label{sec:ShapeDerivatives}
\input{ShapeDerivatives.tex}

\section{Numerical implementation}
\label{sec:SolverAndDiscretization}
\input{GradientDescendMethod_and_Discretization.tex}

\section{Numerical results}
\label{sec:NumericalResults}
\input{NumericalResults.tex}

\section{Conclusion}
\label{sec:Conclusion}
\input{Conclusion.tex}

\section*{References}
\bibliographystyle{elsarticle-harv}
\bibliography{newBibliography.bib}{}

\end{document}

%% file: ALL_ABBR.tex
\newcommand{\abbrStartingCharacteristics}{-}
\newcommand{\abbrEndingCharacteristics}{+}

\newcommand{\surfaceCompDomainStartingCharacteristics}{\denoteSurfacePart{\surfaceCompDomain}{\abbrStartingCharacteristics}}
\newcommand{\surfaceCompDomainEndingCharacteristics}{\denoteSurfacePart{\surfaceCompDomain}{\abbrEndingCharacteristics}}

  \newcommand{\meanCurvature}{h}

\newcommand{\willmoreEnergy}{W}  
\newcommand{\ErosionSquaredIntegral}{E}

\newcommand{\spaceVarTwo}{\vect{y}}

\newcommand{\LTwoDomain}{\lebesqueIntegrableFunctFullDef{2}{\domain}}

\newcommand{\LTwoBoundary}{\lebesqueIntegrableFunctFullDef{2}{\surfaceCompDomain}}

\newcommand{\multiindexDerivative}{\mu}

  \newcommand{\pertubationLetter}{\theta}

  \newcommand{\lipschitzContinuousFunctionsOneRThree}{{C^{1,1}(\R^{3})}^3}

  \newcommand{\lipschitzContinuousFunctionsTwoRThree}{{C^{2,1}(\R^{3})}^3}
  \newcommand{\lipschitzContinuousFunctionsGeneralRThree}{{C^{\orderShapeDerivative,1}(\R^{3})}^3}
    \newcommand{\orderShapeDerivative}{k}
  \newcommand{\normLipschitzContinuousFunctionsTwoRThree}[1]{\norm{#1}_{\lipschitzContinuousFunctionsTwoRThree}}
  \newcommand{\functionSpacePertubationsAdmissableTwoRThree}{\bm{\Theta}^{2}(\R^{3})}

\newcommand{\HOneDomain}{\weakDifferentiableUpToOrder{\domain}
{1}}
  \newcommand{\HOneDomainVect}{\HOneDomain^3}

\newcommand{\HOneHalfBoundary}{\weakDifferentiableUpToOrder{\surfaceCompDomain}{1/2}}
\newcommand{\HOneHalfBoundaryVect}{\weakDifferentiableUpToOrder{\surfaceCompDomain}{1/2}^3}

\newcommand{\shapesAdmissable}{\mathcal{A}}

\newcommand{\trace}[1]{tr(#1)}

\newcommand{\abbrDeformable}{f}

\newcommand{\surfaceCompDomainOutflow}{\denoteSurfacePart{\surfaceCompDomain}{\abbrOutflow}}
\newcommand{\surfaceCompDomainDeformable}{\denoteSurfacePart{\surfaceCompDomain}{\abbrDeformable}}

\newcommand{\IntegralSurfacePartDomainPlain}[2]{\IntegralOne{#1}{#2}{\integrandSurfaceIntegral}}

\newcommand{\tangentialDivergence}[1]{\nabla_{\surfaceDomain}\cdot #1}

\newcommand{\tangentialJacobian}[1]{D_{\surfaceDomain} #1}

\newcommand{\constantRegularizationOne}{c_1}

\newcommand{\costFunctionalLetter}{\integralFunctionalLetter}

\newcommand{\costFunctionalBoundaryFunction}{g}

\newcommand{\jacobian}[1]{\textbf{D}#1}

\newcommand{\genericVectorFunction}{G}

  \newcommand{\gradientInHilbertSpace}{\mathcal{W}} 
    \newcommand{\gradientInHilbertSpaceDiscretized}{\gradientInHilbertSpace^{h}}

\newcommand{\outerNormalTransformed}{\outerNormal_{\domainPertubationLetter}}

\newcommand{\dependenceDomainFromControl}[2]{#1_{#2 \boundaryVelocity}}
  \newcommand{\domainControlled}{\dependenceDomainFromControl{\domain}{\domainPertubationLetter}}

\newcommand{\surfaceDomain}{\Gamma}

\newcommand{\IntegralSurfaceDomainPlain}[1]{\IntegralOne{#1}{\surfaceDomain}{\integrandSurfaceIntegral}}

\newcommand{\normalDerivative}[1]{\partial_{\outerNormalLetter} #1}

\newcommand{\integralFunctionalLetter}{J}

\newcommand{\identityFunction}{I}
    
\newcommand{\boundaryVelocityLetter}{v}
    \newcommand{\boundaryVelocity}{\mathcal{\MakeUppercase{\boundaryVelocityLetter}}}

\newcommand{\domainPertubationLetter}{t}
\newcommand{\transformationFromReferenceDomainLetter}{T}
  \newcommand{\transformationFromReferenceDomain}{\transformationFromReferenceDomainLetter_{\domainPertubationLetter\boundaryVelocity}}

\newcommand{\denoteEulerianDerivativeFunctional}[3]{\dependence{d#1}{#2}[#3]}  
  \newcommand{\denoteEulerianDerivativeFunctionalDomain}[2]{\denoteEulerianDerivativeFunctional{#1}{\domain}{#2}}

\newcommand{\abbrFluid}{f}
\newcommand{\abbrSolid}{p}

\newcommand{\volumePercentageTwo}{\alpha}
\newcommand{\volumePercentageSolid}{\volumePercentageTwo}

\newcommand{\particleIndex}{p}
\newcommand{\denoteParticle}[1]{{#1}_{\particleIndex}}

\newcommand{\diamEq}{\denoteParticle{d}}

\newcommand{\indexNavierStokes}{i}

\newcommand{\velocitySolidLetter}{v}
  \newcommand{\velocitySolid}{\vect{\velocitySolidLetter}}

  \newcommand{\velocitySolidComponent}{\vectorComponent{\velocitySolidLetter}{\indexNavierStokes}}

\newcommand{\velocityFluidLetter}{u}
  \newcommand{\velocityFluid}{\vect{\velocityFluidLetter}}
  \newcommand{\velocityFluidComponent}{\vectorComponent{\velocityFluidLetter}{\indexNavierStokes}}

\newcommand{\densityLetter}{\rho}
  \newcommand{\densityFluid}{\denoteComponent{\densityLetter}{\abbrFluid}}
  \newcommand{\densitySolid}{\denoteComponent{\densityLetter}{\abbrSolid}}

\newcommand{\viscosityLetter}{\mu}
  \newcommand{\viscosityFluid}{\denoteComponent{\viscosityLetter}{\abbrFluid}}
 
\newcommand{\gravityLetter}{g}
  \newcommand{\gravity}{\vect{\gravityLetter}}

\newcommand{\pressureLetter}{p}
  \newcommand{\pressure}{\pressureLetter}

\newcommand{\SchillerNaumannDragTerm}{d_{SN}}
  \newcommand{\weakContributionSchillerNaumannDragTerm}{\tilde{d}_{SN}}
  
  \newcommand{\velocityFluidInflow}{\velocityFluid_{\abbrInflow}}
  \newcommand{\velocitySolidInflow}{\velocitySolid_{\abbrInflow}}
  \newcommand{\volumePercentageSolidInflow}{\volumePercentageSolid_{\abbrInflow}}

\newcommand{\refValue}[1]{{#1}_{ref}}
	\newcommand{\refLength}{\refValue{L}}

	\newcommand{\refVelocity}{\refValue{\velocityFluidLetter}}

\newcommand{\vect}[1]{\boldsymbol{#1}}

\newcommand{\vectorComponent}[2]{#1_{#2}}
\newcommand{\transpose}[1]{#1^{T}}
\newcommand{\inverse}[1]{#1^{-1}}
\newcommand{\inverseTranspose}[1]{#1^{-T}}

\newcommand{\denoteComponent}[2]{{{#1}_{#2}}}
\newcommand{\dependence}[2]{#1(#2)}
\newcommand{\image}[2]{#1(#2)}

\newcommand{\newexpression}[1]{\textbf{#1}}

\newcommand{\reciprocal}[1]{{#1}^{-1}}

 \newcommand{\evaluateAtTwo}[2]{\left.{#1}\right\vert_{#2}}

\newcommand{\derivative}[1]{{#1}'}

	\newcommand{\totalDNoFrac}[2]{\frac{d}{d #2} #1}
	\newcommand{\totalDTwo}[2]{D_{#1} #2}
\newcommand{\multiindexD}[2]{D^{#2} #1}
\newcommand{\partialD}[2]{\partial_{#2} #1}

\newcommand{\directionalDHilbertSpaces}[3]{\totalDTwo{#1}{#2}[#3]}

\newcommand{\divergence}[1]{\nabla\cdot #1}

\newcommand{\gradient}[1]{\nabla #1}

\newcommand{\laplace}[1]{\Delta #1}

\newcommand{\IntegralOne}[3]{\int_{#2} {#1} \mathrm{d} #3}

\newcommand{\integrandSurfaceIntegral}{s}

\newcommand{\defineFunctionLong}[5]{#1\colon #2 &\to     #3\\
				            #4 &\mapsto #5}

\newcommand{\defineFunctionShortSpacesOnly}[3]{#1\colon #2 \to     #3}

\newcommand{\genericFunction}{f}

\newcommand{\N}{\mathbb{N}}

\newcommand{\R}{\mathbb{R}}

	\newcommand{\Rthree}{\R^{3}}

\newcommand{\powerset}[1]{\mathcal{P}(#1)}

\newcommand{\letterLebesqueIntegrableFunctSpace}{L}
	\newcommand{\lebesqueIntegrableFunctFullDef}[2]{\letterLebesqueIntegrableFunctSpace^{#1}(#2)}
	\newcommand{\lebesqueIntegrableAbbrevation}{\letterLebesqueIntegrableFunctSpace^{2}}

\newcommand{\letterWeakDifferentiableFunctSpace}{H}
	\newcommand{\weakDifferentiableUpToOrder}[2]{\letterWeakDifferentiableFunctSpace^{#2}(#1)}

\newcommand{\identityMatrix}{\mathbb{I}}

\newcommand{\abs}[1]{| #1 |}
\newcommand{\norm}[1]{|| #1 ||}
  \newcommand{\normEuclidean}[1]{\norm{#1}_2}

  \newcommand{\normLInfty}[1]{|| #1 ||_{\infty}}

\newcommand{\timeVar}{t}

\newcommand{\spaceVarLetter}{x}

\newcommand{\spaceVar}{\vect{\spaceVarLetter}}

\newcommand{\outerNormalLetter}{n}
  \newcommand{\outerNormal}{\vect{\outerNormalLetter}}

\newcommand{\domain}{\Omega}

\newcommand{\euclideanBasisVector}[1]{e_{#1}}

\newcommand{\closedInterval}[2]{[#1,#2]}

\newcommand{\reynoldsNumber}{Re}
  \newcommand{\reynoldsNumberParticle}{\denoteParticle{\reynoldsNumber}}
  
\newcommand{\froudeNumber}{Fr}
\newcommand{\stokesNumber}{Stk}
	\newcommand{\particleTime}{\timeVar_{\abbrSolid}}

\newcommand{\pecletNumber}{Pe}

\newcommand{\smoothingVelocitySolidNumber}{K}

\newcommand{\functionAngleDependenceImpact}{\zeta}

\newcommand{\functionErosionAtNinetyDegree}{E_{90}}
  \newcommand{\exponentVelocityDependencyOka}{m}
  \newcommand{\exponentDuctileOkaModel}{n_2}
  \newcommand{\exponentBrittleOkaModel}{n_1}
  \newcommand{\hardnessOkaModel}{H_{v}}

  \newcommand{\exponentBrittleOkaModelDuctileCaseVal}{0.78}
  \newcommand{\exponentDuctileOkaModelDuctileCaseVal}{1.25}

\newcommand{\functionalErosionLetter}{e}

\newcommand{\angleErosion}{\gamma}

\newcommand{\depositionFraction}{\eta}

\newcommand{\deanNumber}{De}
\newcommand{\dTube}{d_t}

\newcommand{\rCurvatureBend}{r_b}
\newcommand{\rCurvatureRatio}{R_0}
\newcommand{\velocityFluidMeanInflow}{\velocityFluidLetter_0}

  \newcommand{\abbrInflow}{in}
  \newcommand{\abbrOutflow}{out}

  \newcommand{\abbrWall}{w}

\newcommand{\surfaceCompDomain}{\Gamma}
  \newcommand{\denoteSurfacePart}[2]{#1^{#2}}
  \newcommand{\surfaceCompDomainInflow}{\denoteSurfacePart{\surfaceCompDomain}{\abbrInflow}}

  \newcommand{\surfaceCompDomainWall}{\denoteSurfacePart{\surfaceCompDomain}{\abbrWall}}

\newcommand{\IntegralDomain}[1]{\IntegralOne{#1}{\domain}{\spaceVar}}

\newcommand{\spaceX}{X(\domain)}
\newcommand{\spaceXZero}{X_{0}(\domain)}
\newcommand{\spaceY}{Y(\domain)}
\newcommand{\spaceYZero}{Y_{0}(\domain)}
\newcommand{\spaceM}{M(\domain)}
\newcommand{\spaceZ}{Z(\domain)}
\newcommand{\spaceZZero}{Z_{0}(\domain)}

\newcommand{\spaceBOne}{B_1(\domain)}
\newcommand{\spaceBTwo}{B_2(\domain)}

\newcommand{\spaceE}{E(\domain)}
\newcommand{\spaceF}{F(\domain)}

\newcommand{\spaceEDeformed}{E(\domainControlled)}
\newcommand{\spaceFDeformed}{F(\domainControlled)}

\newcommand{\genericFunctionalArgumentOne}{\kappa}
\newcommand{\genericFunctionalArgumentTwo}{\lambda}
\newcommand{\genericFunctionalArgumentThree}{\mu}

\newcommand{\functionalA}{\textup{a}}
\newcommand{\functionalB}{\textup{b}}
\newcommand{\functionalC}{\textup{c}}

\newcommand{\functionalF}{\textup{f}}
\newcommand{\functionalG}{\textup{g}}

\newcommand{\functionalHVelocityFluidInflow}{\textup{h}_1}
\newcommand{\functionalHVelocitySolidInflow}{\textup{h}_2}
\newcommand{\functionalHVolumePercentageSolidInflow}{\textup{h}_3}

\newcommand{\functionalD}{\textup{d}}

\newcommand{\functionalEq}{\textup{Q}}

\newcommand{\lagrangianSturmApproach}{G}
  \newcommand{\lagrangianSturmApproachTransformed}{\mathscr{\lagrangianSturmApproach}}

\newcommand{\ansatzSpaceSturmApproachLetter}{\Phi}
  \newcommand{\ansatzSpaceSturmApproachVelocity}{\vect{\ansatzSpaceSturmApproachLetter}_{\velocityFluidLetter}}
  \newcommand{\ansatzSpaceSturmApproachPressure}{\ansatzSpaceSturmApproachLetter_{\pressure}}

  \newcommand{\ansatzSpaceSturmApproachVelocitySolid}{\vect{\ansatzSpaceSturmApproachLetter}_{\velocitySolidLetter}}
    \newcommand{\ansatzSpaceSturmApproachVolumePercentageSolid}{\ansatzSpaceSturmApproachLetter_{\volumePercentageSolid}}

    \newcommand{\ansatzSpaceSturmApproachVelocityComponent}{\ansatzSpaceSturmApproachLetter_{\velocityFluidComponent}}
  \newcommand{\ansatzSpaceSturmApproachVelocitySolidComponent}{\ansatzSpaceSturmApproachLetter_{\velocitySolidComponent}}

\newcommand{\testSpaceSturmApproachLetter}{\Psi}
  \newcommand{\testSpaceSturmApproachVelocity}{\vect{\testSpaceSturmApproachLetter}_{\velocityFluidLetter}}
  \newcommand{\testSpaceSturmApproachPressure}{\testSpaceSturmApproachLetter_{\pressure}}

  \newcommand{\testSpaceSturmApproachVelocitySolid}{\vect{\testSpaceSturmApproachLetter}_{\velocitySolidLetter}}
  \newcommand{\testSpaceSturmApproachVolumePercentageSolid}{\testSpaceSturmApproachLetter_{\volumePercentageSolid}}
  \newcommand{\testSpaceSturmApproachVelocityFluidInflowCondition}{\vect{\testSpaceSturmApproachLetter}_{{\velocityFluidLetter}_{in}}}
  \newcommand{\testSpaceSturmApproachVelocitySolidInflowCondition}{\vect{\testSpaceSturmApproachLetter}_{{\velocitySolidLetter}_{in}}}
  \newcommand{\testSpaceSturmApproachVolumePercentageSolidInflowCondition}{\testSpaceSturmApproachLetter_{{\volumePercentageSolid}_{in}}}
  
\newcommand{\solutionSturmApproachLetter}{q}
  \newcommand{\solutionSturmApproachDomain}{\solutionSturmApproachLetter}
  \newcommand{\solutionSturmApproachVelocityFluid}{\vect{\solutionSturmApproachDomain}_{\velocityFluidLetter}}
  \newcommand{\solutionSturmApproachPressure}{\solutionSturmApproachDomain_{\pressure}}

  \newcommand{\solutionSturmApproachVelocitySolid}{\vect{\solutionSturmApproachDomain}_{\velocitySolidLetter}}
  \newcommand{\solutionSturmApproachVolumePercentageSolid}{\solutionSturmApproachDomain_{\volumePercentageSolid}}
  
  \newcommand{\solutionSturmApproachVelocityFluidComponent}{\solutionSturmApproachDomain_{\velocityFluidComponent}}
  \newcommand{\solutionSturmApproachVelocitySolidComponent}{\solutionSturmApproachDomain_{\velocitySolidComponent}}

\newcommand{\adjointSolutionSturmApproachLetter}{\lambda}
  \newcommand{\adjointSolutionSturmApproachDomain}{\adjointSolutionSturmApproachLetter}
  \newcommand{\adjointSolutionSturmApproachVelocityFluid}{\vect{\adjointSolutionSturmApproachDomain}_{\velocityFluidLetter}}
  \newcommand{\adjointSolutionSturmApproachPressure}{\adjointSolutionSturmApproachDomain_{\pressure}}

    \newcommand{\adjointSolutionSturmApproachVelocitySolid}{\vect{\adjointSolutionSturmApproachLetter}_{\velocitySolidLetter}}
    \newcommand{\adjointSolutionSturmApproachVolumePercentageSolid}{\adjointSolutionSturmApproachLetter_{\volumePercentageSolid}}
    \newcommand{\adjointSolutionSturmApproachVelocityFluidInflowCondition}{\vect{\adjointSolutionSturmApproachLetter}_{{\velocityFluidLetter}_{in}}}
    \newcommand{\adjointSolutionSturmApproachVelocitySolidInflowCondition}{\vect{\adjointSolutionSturmApproachLetter}_{{\velocitySolidLetter}_{in}}}
      \newcommand{\adjointSolutionSturmApproachVolumePercentageSolidInflowCondition}{\adjointSolutionSturmApproachLetter_{{\volumePercentageSolid}_{in}}}

\newcommand{\determinantJacobianLetter}{\xi}
  \newcommand{\detTrafo}{\determinantJacobianLetter}
    \newcommand{\detTrafoEvaluated}{\dependence{\detTrafo}{\domainPertubationLetter}}

  \newcommand{\detTrafoBdry}{\determinantJacobianLetter_\surfaceCompDomain}
    \newcommand{\detTrafoBdryEvaluated}{\dependence{\detTrafoBdry}{\domainPertubationLetter}}

    \newcommand{\detTrafoBdryPrime}{\derivative{\detTrafoBdry}}
    
     \newcommand{\detTrafoBdryPrimeAtZero}{\dependence{\detTrafoBdryPrime}{0}}
    
  \newcommand{\jacTrafo}{\jacobian{\transformationFromReferenceDomain}}
  \newcommand{\invTranspJacTrafo}{M}
    \newcommand{\invTranspJacTrafoEvaluated}{\dependence{\invTranspJacTrafo}{\domainPertubationLetter}}

  \newcommand{\invJacTrafoEvaluated}{\dependence{\transpose{\invTranspJacTrafo}}{\domainPertubationLetter}}

  \newcommand{\jacDefVelo}{\jacobian{\boundaryVelocity}}

 \newcommand{\spaceForGradientOfComposedScalarFunction}{W^{1,1}_{loc}(\R^{3})}
\newcommand{\spaceForGradientOfComposedVectorFunction}{{W^{1,1}_{loc}(\R^{3})}^3}
 
 \newcommand{\hilbertSpaceDeformation}{\mathcal{D}}
  
  \newcommand{\hilbertSpaceDeformationFiniteDimensional}{\R^{3\expDiscrHilbertSpaceDef}}
    \newcommand{\expDiscrHilbertSpaceDef}{k}

  \newcommand{\bilinearFormProjection}{\textup{A}}
  
  \newcommand{\strainTensor}{\sigma}
  \newcommand{\stressTensor}{\epsilon}
  \newcommand{\lameParameterOne}{\lambda}
  \newcommand{\lameParameterTwo}{\mu}
    \newcommand{\lameParameterTwoPDE}{{\lameParameterTwo}^{*}}
    \newcommand{\lameParameterTwoPDEMin}{\lameParameterTwoPDE_{min}}
    \newcommand{\lameParameterTwoPDEMax}{\lameParameterTwoPDE_{max}}
  \newcommand{\PoissonRation}{\nu}

 \newcommand{\normalForceOperator}{N}
  
 \newcommand{\boundaryVelocityCorrection}{\Pi}
  \newcommand{\boundaryVelocityCorrectionDiscretized}{\boundaryVelocityCorrection^{h}}
 \newcommand{\normalForceLagrangeParameter}{F}
    \newcommand{\normalForceLagrangeParameterDiscretized}{\normalForceLagrangeParameter^{h}}
 \newcommand{\normalForceLagrangeParameterTestFunction}{E}
 \newcommand{\gradientInHilbertSpaceCorrected}{{G}}
  \newcommand{\gradientInHilbertSpaceCorrectedDiscretized}{\gradientInHilbertSpaceCorrected^{h}}

  \newcommand{\lbfgsDescentDirectionDiscrete}{q^{h}}
 
 \newcommand{\stepsizeBacktracking}{t}

 \newcommand{\iterBFGS}{j}
 \newcommand{\denoteIteration}[1]{#1_\iterBFGS}
 
 \newcommand{\denoteNextIteration}[1]{#1_{\iterBFGS+1}}

 \newcommand{\meshVertices}{X^{h}}


%% file: Abstract.tex
In this paper we consider a shape optimization problem for the minimization of the erosion,
that is caused by the impact of inert particles onto the walls of a bended pipe.
Using the continuous adjoint approach, 
we formally compute the shape derivative of the optimization problem, 
which is based on a one-way coupled, fully Eulerian description of a monodisperse particle jet, that is transported in a carrier fluid.
We validate our approach by numerically optimizing a three-dimensional pipe segment with respect to a single particle species using a gradient descent method,
and show, that the 
erosion rates on the optimized geometry are reduced with respect to the initial bend for a broader range of particle Stokes numbers.

%% file: Introduction.tex
The erosion caused by a dilute suspension of inert particles in a carrier fluid due to particle-wall collisions is a lifetime determining factor especially for bended parts of pipe systems \cite{mills2003pneumatic}.
Its reduction through modifications of the geometry is therefore an active field of research, 
see e.g. \cite{dos2016reducing,duarte2017innovative,fan2002antierosion,zhu2018numerical} and the references therein.
The prevalent numerical method for such studies is the
Lagrangian particle-tracking method, where the erosion rate is determined from the velocities and impact angles of a discrete number of representative particles.

While effects such as rebound and secondary impact are not easily incorporated in an
Eulerian treatment of the particulate phase \cite{sachdev2007numerical, saurel1994two, slater2001calculation},
this approach avoids the need to choose an adequate sample size
due to the presence of continuously varying locally averaged particle variables.
Additionally, the minimization of the erosion rate can in this setting 
be formulated as a PDE-constrained shape optimization problem \cite{hinze2008optimization,troltzsch2005optimale},
so that continuous adjoint techniques can be applied to improve the initial geometry. 

To make an advance in this field, 
we consider in this work a shape optimization problem for a bended pipe,
where the erosion due to the primary impact of a dilute monodisperse particle jet transported in laminar air flow is to be minimized.
In order to do so, we only allow for deformations of the bend without introducing additional replaceable parts such as twisted baffles upstream of the curved segment, 
as done e.g. in \cite{dos2016reducing}.
We describe the particle transport with the one-way coupled Eulerian model from \cite{bourgault2000three,bourgault1999finite}
and formally compute the shape derivative of a cost functional 
based on a volume-averaged formulation of the erosion rate 
predicted by the OKA model \cite{oka2005practical}.
To the best of the authors' knowledge, this provides a new result, 
which can easily be extended to other commonly encountered erosion models such as the ones described in \cite{finnie1972,zhang2007comparison}, 
since only the partial derivatives of the erosion rate with respect to the state variables have to be substituted.
In order to obtain a descent direction, the shape derivative is projected 
based on the equations of linear elasticity and 
used within a gradient descent method for the numerical optimization of the initial bend.

This paper is organized in the following way:
In \cref{sec:ProblemFormulation} we introduce the Eulerian particle model and the optimal shape design problem.
\Cref{sec:ShapeDerivatives} is devoted to the notation from shape calculus, 
with which we compute the shape derivative of the cost functional in \Cref{thm:ShapeDerivative}.
Our approach to obtain smooth volume mesh deformations from a Riesz projection of the shape derivative,
as well as the gradient descent method \Cref{alg:GradientDescentMethod},
are described in \cref{sec:SolverAndDiscretization}.
In this section we also comment on the discretization 
of the PDEs with finite elements, the stabilization of the discrete formulation, 
and the numerical solution procedure.
In \cref{sec:NumericalResults} we validate our approach by comparing the calculated impact rates of the Eulerian particle model for various Stokes numbers with reference values from the literature.
We then optimize the three-dimensional reference geometry with respect to the erosion rate caused by a selected particle species with an intermediate Stokes number
and show,
that the optimized geometry experiences less erosion than the initial bend for a broader range of particle species.

%% file: ProblemFormulation.tex
In this section we introduce the
fluid equations and the 
Eulerian particle transport model, 
which act as PDE-constraints for the 
shape optimization problem introduced in \cref{sec:OptimizationProblem}.
The model equations 
together with a description of the geometry of the test case considered in \cref{sec:NumericalResults}
are given in \cref{sec:MathematicalModel}.
Since we use a finite element approach for the discretization of the PDEs,
we derive the weak formulation of the problem in \cref{sec:WellPosednessWeakFormulation}.

\subsection{Mathematical model}
\label{sec:MathematicalModel}

{ 
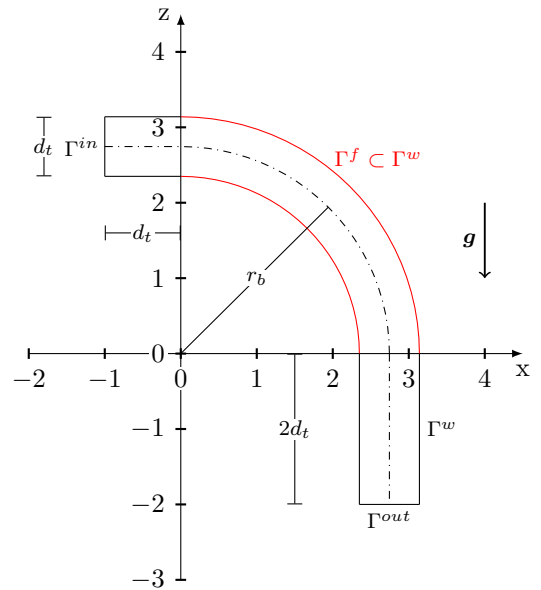
\begin{wrapfigure}{r}{0.45\textwidth}
  \begin{center}
    \input{Pictures/tikzpicture_bend.tex}
  \end{center}
  \caption{Two-dimensional description of the bended pipe with deformable boundary}
  \label{fig:PlotBendedPipe}
\end{wrapfigure}

Consider the flow of a fluid containing 
small, spherical, solid particles
in a bended pipe $\domain\subset\Rthree$ with circular cross section of diameter $\dTube$
and radius of curvature $\rCurvatureBend$ 
as in \Cref{fig:PlotBendedPipe}.
The boundary $\surfaceCompDomain$ of $\domain$ is subdivided in the form
$\surfaceCompDomain = \surfaceCompDomainInflow \dot\cup \surfaceCompDomainWall \dot\cup \surfaceCompDomainOutflow$,
where the individual parts are assumed to be of class $C^2$.
Our goal is to optimize the bended segment 
$\surfaceCompDomainDeformable \subset \surfaceCompDomainWall$ 
connecting the inlet and outlet, 
so that for prescribed mass flow rates at  $\surfaceCompDomainInflow$,
the predicted erosion rate at $\surfaceCompDomainWall$ can be reduced.

By assuming a low particle load, the fluid and particle equations can be decoupled, 
since the presence of the particles does in this case not considerably influence the carrier fluid.
Given a reference velocity $\refVelocity$ of the flow, the reference length $\refLength:=\dTube$, 
as well as the density $\densityFluid$ and viscosity $\viscosityFluid$ of the fluid 
and gravity $\gravityLetter$ acting in direction $\gravity = - \euclideanBasisVector{z}$,
we define the Reynolds number 
$\reynoldsNumber:=\densityFluid \refVelocity \refLength / \viscosityFluid$ 
and the Froude number 
$\froudeNumber:=\refVelocity/\sqrt{\refLength \gravityLetter}$.
Additionally we introduce the 
Dean number 
$\deanNumber = \reynoldsNumber/\sqrt{\rCurvatureRatio}$
with the curvature ratio $\rCurvatureRatio=2\,\rCurvatureBend/\dTube$,
an important dimensionless number for the characterization of fluid flow in curved pipes \cite{dean1927}.

The fluid velocity $\velocityFluid\in\Rthree$ and pressure $\pressure \in \R$
are described with the 
non-dimensionalized, stationary, incompressible, isothermal Navier-Stokes equations
\begin{subequations}
\begin{align}
  \label{eq:NavierStokesEquations}
   \transpose{{\velocityFluid}} \transpose{\gradient{{\velocityFluid}}} 
  +\gradient{{\pressure}} 
  - \reciprocal{\reynoldsNumber} \laplace{{\velocityFluid}}
  &=  {\froudeNumber}^{-2} \gravity \textup{ in } \domain, \\
  \label{eq:NavierStokesIncompressibility}
  \divergence{\velocityFluid} &= 0 \textup{ in } \domain .
\end{align}
\end{subequations}
We impose the boundary conditions 
\begin{subequations}
\begin{align}
    \label{eq:BoundaryConditionNavierStokesInflow}
      \velocityFluid &= \velocityFluidInflow \textup{ on } \surfaceCompDomainInflow, \\
    \label{eq:BoundaryConditionNavierStokesWall}
      \velocityFluid &= 0 \textup{ on } \surfaceCompDomainWall, \\
    \label{eq:BoundaryConditionNavierStokesOutflow}
      \pressure \outerNormal - \reciprocal{\reynoldsNumber}\normalDerivative{\velocityFluid} 
      &= 0  \textup{ on } \surfaceCompDomainOutflow, 
\end{align}
\end{subequations}
where $\velocityFluidInflow$ is the parabolic inflow profile of a fully-developed laminar pipe flow, 
and the reference velocity is chosen to be 
$\refVelocity = 0.5 \max{\normEuclidean{\velocityFluidInflow}}$.
} 

Important dimensionless numbers characterizing the sensitivity of 
the particles  of diameter $\diamEq$ and density $\densitySolid$ to the local flow pattern are 
the 
particle response time 
$\particleTime:= \densitySolid  \diamEq^2/(18 \viscosityFluid)$
and the 
Stokes number 
$\stokesNumber:= \particleTime\refVelocity / (0.5\, \refLength)$.
Particles with $\stokesNumber \ll 1$ closely follow the fluid streamlines, 
such that almost no particle-wall interaction takes place.
For increasing values of $\stokesNumber$, the fluid drag exerted onto the particles becomes less dominant in comparison to the inertial forces,
so that the particle trajectories start to deviate from the fluid streamlines. 
In this case a certain amount of particles 
will come in contact with the wall $\surfaceCompDomainWall$ and cause erosion upon impact.

In this paper, we describe the particle phase with one locally averaged velocity field
following \cite{bourgault1999finite, bourgault2000three, honsek2008eulerian, slater2001calculation}.
In this model, 
the volume-averaged, dimensionless particle velocity $\velocitySolid$ is obtained from
\begin{equation}
 \label{eq:ParticleEquationVelocity}
\transpose{{\velocitySolid}} \transpose{\gradient{{\velocitySolid}}} 
  + 2\,\reciprocal{\stokesNumber}\dependence{\SchillerNaumannDragTerm}{\velocityFluid,\velocitySolid}({\velocitySolid}-{\velocityFluid})
  - \reciprocal{\smoothingVelocitySolidNumber} \laplace{\velocitySolid}
  = {\froudeNumber}^{-2} \gravity \textup{ in } \domain,
\end{equation}	
where the second order term with $\smoothingVelocitySolidNumber\gg 1$ is added 
for regularity reasons only.
The particle Reynolds number
$
\dependence{\reynoldsNumberParticle}{\velocityFluid,\velocitySolid}:=
\densityFluid \normEuclidean{\velocitySolid-\velocityFluid} \diamEq/\viscosityFluid
$
is assumed to be below $1000$ in the entire domain,
such that the fluid drag exerted on the particulate phase 
can be modeled with the Schiller-Naumann correlation \cite{schiller1935}
\begin{equation*}
 \dependence{\SchillerNaumannDragTerm}{\velocityFluid,\velocitySolid} :=
 1+0.15 \dependence{\reynoldsNumberParticle}{\velocityFluid,\velocitySolid}^{0.687}.
\end{equation*}
\Cref{eq:ParticleEquationVelocity} is completed with the boundary conditions \cite{bourgault1999finite}
\begin{subequations}
\begin{align}
\label{eq:BdryConditionsSolidVelocityInflow}
 \velocitySolid &= \velocitySolidInflow \textup{ on } \surfaceCompDomainInflow, \\
\label{eq:BdryConditionsSolidVelocityNonInflow}
  \normalDerivative{\velocitySolid} &=0 \textup{ on } \surfaceCompDomain \setminus \surfaceCompDomainInflow ,
\end{align}
\end{subequations}
with the particle velocity profile $\velocitySolidInflow$,
where we choose $\velocityFluidInflow=\velocitySolidInflow$ for simplicity.

Given a solution of \cref{eq:ParticleEquationVelocity,eq:BdryConditionsSolidVelocityInflow,eq:BdryConditionsSolidVelocityNonInflow}, 
the volume percentage occupied by the particles is obtained from the advection-diffusion equation
\begin{equation}
  \label{eq:ParticleEquationTransport}
 \divergence{({\volumePercentageSolid}
  {\velocitySolid} )}
  -\reciprocal{\pecletNumber} \laplace{{\volumePercentageSolid}}
  = 0 \textup{ in } \domain,
\end{equation}
where similarly to \cref{eq:ParticleEquationVelocity}, 
artificial diffusion depending on the Peclet number $\pecletNumber \gg 1$ is added
for regularization only.
Owing to the theory of characteristics, 
the boundary conditions for $\volumePercentageSolid$ have to depend on the sign of $\transpose{\velocitySolid}\outerNormal$.
Following \cite{bourgault1999finite},
we set
\begin{equation*}
\dependence{\surfaceCompDomainStartingCharacteristics}{\velocitySolid} 
	 := \{ \spaceVar \in \surfaceCompDomain \, \vert \, \transpose{\dependence{\velocitySolid}{\spaceVar}} \dependence{\outerNormal}{\spaceVar} \le 0
	 \}, \quad
\dependence{\surfaceCompDomainEndingCharacteristics}{\velocitySolid} := \surfaceCompDomain \setminus \dependence{\surfaceCompDomainStartingCharacteristics}{\velocitySolid},
\end{equation*}
and consider \cref{eq:ParticleEquationTransport} together with 
\begin{subequations}
\begin{align}
\label{eq:BdryConditionsSolidVolumePercentageInflow}
   \volumePercentageSolid 
  &= \volumePercentageSolidInflow \textup{ on } \surfaceCompDomainInflow, \\
\label{eq:BdryConditionsSolidVolumePercentageStartingCharacteristics}
\volumePercentageSolid 
&= 0  \textup{ on }
\dependence{\surfaceCompDomainStartingCharacteristics}{\velocitySolid}\setminus\surfaceCompDomainInflow, \\
\label{eq:BdryConditionsSolidVolumePercentageEndingCharacteristics}
  \normalDerivative{\volumePercentageSolid} &= 0 \textup{ on } \dependence{\surfaceCompDomainEndingCharacteristics}{\velocitySolid} .
\end{align}
\end{subequations}
Here $\volumePercentageSolidInflow>0$ is the distribution of particles at the inlet,
which is taken to be constant throughout this paper.
The Dirichlet boundary condition \cref{eq:BdryConditionsSolidVolumePercentageStartingCharacteristics}
avoids unphysical inflow from so-called \newexpression{shadowed zones}, 
where particles do not make contact with the boundary \cite{slater2001calculation}
and the Neumann condition \cref{eq:BdryConditionsSolidVolumePercentageEndingCharacteristics} 
leads to a flux, 
\begin{equation*}
\IntegralSurfacePartDomainPlain{\volumePercentageSolid (\transpose{\velocitySolid}  \outerNormal)}{\surfaceCompDomainWall \cap \dependence{\surfaceCompDomainEndingCharacteristics}{\velocitySolid}},
\end{equation*}
which can be used to identified the particle impact rate.

In contrast to Lagrangian methods, 
particle rebound and crossing particle trajectories pose severe difficulties in case only one particle velocity  is used in the Eulerian model \cite{slater2001calculation}.
To circumvent this restriction,
it has been proposed \cite{sachdev2007numerical, saurel1994two} 
to consider several particle species for different directions
and more recent work in this field focuses on multi-velocity formulations obtained from kinetic theory \cite{desjardins2008quadrature,forgues2019higher,fox2009higher}.
However, due to their complexity, 
these models do not lend themselves easily to the continuous adjoint approach of optimization
and since the focus of this work lies 
on the shape optimization of a flow situation with one dominating convective transport direction, 
we chose to use the simpler model governed by \cref{eq:ParticleEquationVelocity,eq:ParticleEquationTransport} 
subject to \cref{eq:BdryConditionsSolidVelocityInflow,eq:BdryConditionsSolidVelocityNonInflow,eq:BdryConditionsSolidVelocityInflow,eq:BdryConditionsSolidVolumePercentageInflow,eq:BdryConditionsSolidVolumePercentageStartingCharacteristics,eq:BdryConditionsSolidVolumePercentageEndingCharacteristics}, 
while being aware of its limitations. 

\subsection{Weak formulation}
\label{sec:WellPosednessWeakFormulation}

The weak formulation of the Eulerian model introduced in \cref{sec:MathematicalModel}
is based on the function spaces 
\begin{align*}
 \spaceX &:= \{ \ansatzSpaceSturmApproachVelocity \in \HOneDomainVect \,\vert\, \evaluateAtTwo{\ansatzSpaceSturmApproachVelocity}{\surfaceCompDomainWall} = 0\}, \quad
 \spaceXZero := \{ \testSpaceSturmApproachVelocity \in \spaceX \,\vert\, 
  \evaluateAtTwo{\testSpaceSturmApproachVelocity }{\surfaceCompDomainInflow} = 0 \}, \\
  \spaceM &:=\LTwoDomain, \quad
  \spaceY := \HOneDomainVect, \quad
  \spaceYZero := \{ \testSpaceSturmApproachVelocitySolid \in \spaceY \,\vert\,
    \evaluateAtTwo{\testSpaceSturmApproachVelocitySolid}{\surfaceCompDomainInflow} = 0\}, \\
  \spaceZ &:= \{\ansatzSpaceSturmApproachVolumePercentageSolid \in \HOneDomain \vert
	      \evaluateAtTwo{\ansatzSpaceSturmApproachVolumePercentageSolid}{\dependence{\surfaceCompDomainStartingCharacteristics}{\solutionSturmApproachVelocitySolid} \setminus \surfaceCompDomainInflow} = 0\}, \quad
  \spaceZZero := \{ \testSpaceSturmApproachVolumePercentageSolid \in \spaceZ \,\vert\, 
      \evaluateAtTwo{\testSpaceSturmApproachVolumePercentageSolid}{\surfaceCompDomainInflow}=0\}, \\
  \spaceBOne &:= \HOneHalfBoundaryVect, \quad 
  \spaceBTwo := \HOneHalfBoundary,
\end{align*}
equipped with the standard  $\lebesqueIntegrableAbbrevation$-norms,
where the restriction of functions onto boundary parts is meant in the sense of traces.
Since the boundary conditions \cref{eq:BdryConditionsSolidVolumePercentageStartingCharacteristics,eq:BdryConditionsSolidVolumePercentageEndingCharacteristics} for the volume percentage depend on the sign of the normal component of the particle velocity,
the solution $\solutionSturmApproachVelocitySolid\in\spaceY$ of the weak form of \cref{eq:ParticleEquationVelocity,eq:BdryConditionsSolidVelocityInflow,eq:BdryConditionsSolidVelocityNonInflow},
which will be given shortly,
is already required for the definition of $\spaceZ$ and $\spaceZZero$.

For the sake of brevity we further denote by 
\begin{align*}
  \spaceE &:=\spaceX \times \spaceM \times \spaceY \times \spaceZ, \\
  \spaceF &:= \spaceXZero \times \spaceM \times \spaceYZero \times \spaceZZero \times \spaceBOne \times \spaceBOne \times \spaceBTwo, 
\end{align*}
the function spaces for the forward and adjoint variables respectively 
and define the mappings 
 \begingroup
\allowdisplaybreaks
\begin{align*}
 \dependence{\functionalA}{\domain,\genericFunctionalArgumentOne,\genericFunctionalArgumentTwo} 
 &:= \IntegralDomain{\gradient{\genericFunctionalArgumentOne} : \gradient{\genericFunctionalArgumentTwo}}, \quad 
 \dependence{\functionalB}{\domain,\genericFunctionalArgumentOne,\genericFunctionalArgumentTwo}
 := \IntegralDomain{\genericFunctionalArgumentTwo (\divergence{\genericFunctionalArgumentOne})}, \\
 \dependence{\functionalC}{\domain, \genericFunctionalArgumentOne, \genericFunctionalArgumentTwo, \genericFunctionalArgumentThree}
 &:= \IntegralDomain{\transpose{\genericFunctionalArgumentOne}  \transpose{\gradient{\genericFunctionalArgumentTwo}}  \genericFunctionalArgumentThree}, \quad
 \dependence{\functionalD}{\domain, \genericFunctionalArgumentOne, \genericFunctionalArgumentThree,\genericFunctionalArgumentTwo}
 := \IntegralDomain{\divergence{(\genericFunctionalArgumentOne \genericFunctionalArgumentThree)\genericFunctionalArgumentTwo}}, \\
\dependence{\functionalF}{\domain,\genericFunctionalArgumentOne,\genericFunctionalArgumentTwo,\genericFunctionalArgumentThree}
 &:= \IntegralDomain{\dependence{\SchillerNaumannDragTerm}{\genericFunctionalArgumentOne,\genericFunctionalArgumentTwo}\transpose{(\genericFunctionalArgumentOne-\genericFunctionalArgumentTwo)}  \genericFunctionalArgumentThree}, \quad
 \dependence{\functionalG}{\domain, \genericFunctionalArgumentThree}
 := \IntegralDomain{\transpose{\gravity}  \genericFunctionalArgumentThree}, \\
 \dependence{\functionalHVelocityFluidInflow}{\domain,\genericFunctionalArgumentOne,\genericFunctionalArgumentTwo}
  &:= \IntegralSurfacePartDomainPlain{\transpose{(\genericFunctionalArgumentOne - \velocityFluidInflow)} \genericFunctionalArgumentTwo}{\surfaceCompDomainInflow}, \quad 
 \dependence{\functionalHVelocitySolidInflow}{\domain,\genericFunctionalArgumentOne,\genericFunctionalArgumentTwo}
 := \IntegralSurfacePartDomainPlain{\transpose{(\genericFunctionalArgumentOne - \velocitySolidInflow)} \genericFunctionalArgumentTwo}{\surfaceCompDomainInflow}, \\
 \dependence{\functionalHVolumePercentageSolidInflow}{\domain,\genericFunctionalArgumentOne,\genericFunctionalArgumentTwo}
 &:= \IntegralSurfacePartDomainPlain{(\genericFunctionalArgumentOne - \volumePercentageSolidInflow) \genericFunctionalArgumentTwo}{\surfaceCompDomainInflow}.
\end{align*}
\endgroup

With these definitions, 
the weak formulation of the stationary Eulerian particle transport model
reads:
\begin{align}
\begin{split}
\textup{Seek }
\solutionSturmApproachLetter = (\solutionSturmApproachVelocityFluid,\solutionSturmApproachPressure,\solutionSturmApproachVelocitySolid, \solutionSturmApproachVolumePercentageSolid)  \in \spaceE
\textup{, such that} & \\
\label{eq:WeakForwardPDESystem}
  \dependence{\functionalC}{\domain,\solutionSturmApproachVelocityFluid,\solutionSturmApproachVelocityFluid,\testSpaceSturmApproachVelocity}
  + \reciprocal{\reynoldsNumber}\dependence{\functionalA}{\domain,\solutionSturmApproachVelocityFluid,\testSpaceSturmApproachVelocity} 
  - \dependence{\functionalB}{\domain,\testSpaceSturmApproachVelocity,\solutionSturmApproachPressure}
 &= {\froudeNumber}^{-2}\dependence{\functionalG}{\domain,\testSpaceSturmApproachVelocity},\\ 
 \dependence{\functionalB}{\domain,\solutionSturmApproachVelocityFluid,\testSpaceSturmApproachPressure}
 &= 0,\\
 \dependence{\functionalHVelocityFluidInflow}{\domain,\solutionSturmApproachVelocityFluid,\testSpaceSturmApproachVelocityFluidInflowCondition}
 &= 0,\\
 \dependence{\functionalC}{\domain,\solutionSturmApproachVelocitySolid,\solutionSturmApproachVelocitySolid,\testSpaceSturmApproachVelocitySolid}
 +\reciprocal{\smoothingVelocitySolidNumber}
   \dependence{\functionalA}{\domain,\solutionSturmApproachVelocitySolid,\testSpaceSturmApproachVelocitySolid} 
   + 2\reciprocal{\stokesNumber} \dependence{\functionalF}{\domain, \solutionSturmApproachVelocityFluid,\solutionSturmApproachVelocitySolid,\testSpaceSturmApproachVelocitySolid}
   &= {\froudeNumber}^{-2}\dependence{\functionalG}{\domain,\testSpaceSturmApproachVelocitySolid}, \\ 
    \dependence{\functionalHVelocitySolidInflow}{\domain,\solutionSturmApproachVelocitySolid,\testSpaceSturmApproachVelocitySolidInflowCondition}
 &= 0, \\
   \dependence{\functionalD}{\domain,\solutionSturmApproachVolumePercentageSolid,\solutionSturmApproachVelocitySolid,\testSpaceSturmApproachVolumePercentageSolid}
   +\reciprocal{\pecletNumber}\dependence{\functionalA}{\solutionSturmApproachVolumePercentageSolid,\testSpaceSturmApproachVolumePercentageSolid}
   &= 0, \\
     \dependence{\functionalHVolumePercentageSolidInflow}{\domain,\solutionSturmApproachVolumePercentageSolid,\testSpaceSturmApproachVolumePercentageSolidInflowCondition}
     &= 0, \\
\textup{holds for all }
(\testSpaceSturmApproachVelocity,\testSpaceSturmApproachPressure,\testSpaceSturmApproachVelocitySolid,\testSpaceSturmApproachVolumePercentageSolid,\testSpaceSturmApproachVelocityFluidInflowCondition,\testSpaceSturmApproachVelocitySolidInflowCondition,\testSpaceSturmApproachVolumePercentageSolidInflowCondition) \in \spaceF. &
     \end{split}
 \end{align}

\subsection{Optimization problem}
\label{sec:OptimizationProblem}

If the particles are treated in a Lagrangian manner, 
the erosion caused by their impact onto the wall is generally predicted with 
either the model of Finnie \cite{finnie1972}, 
the Oka model \cite{oka2005practical} or the E/CRC model \cite{zhang2007comparison}. 
A common trait of these models is, that for a given pipe material and given particle parameters  
the predicted erosion pattern continuously depends on the impact angle and velocity as well as the amount of particles hitting the wall.
In the following, we restrict ourselves to the dimensionless Eulerian formulation 
\begin{equation}
\label{eq:ErosionRateOKA}
 \dependence{\functionalErosionLetter}{\volumePercentageSolid,\velocitySolid,\outerNormal} 
 = \volumePercentageSolid (\velocitySolid \cdot \outerNormal) 
  \dependence{\functionErosionAtNinetyDegree}{\velocitySolid} \dependence{\functionAngleDependenceImpact}{\dependence{\angleErosion}{\velocitySolid,\outerNormal}},
  \quad \textup{ if } \transpose{\velocitySolid}\outerNormal \ge 0, 
\end{equation}
of the erosion rate from \cite{oka2005practical}
and note, that the other erosion models can be expressed similarly.
Here
\begin{equation}
\label{eq:AngleDependentFunctionOka}
 \dependence{\functionErosionAtNinetyDegree}{\velocitySolid} = \normEuclidean{\velocitySolid}^{\exponentVelocityDependencyOka},
 \quad 
 \exponentVelocityDependencyOka > 0,
\end{equation}
is the erosion rate for orthogonal impact with respect to the tangential plane,
and the functions
\begin{align*}
\dependence{\functionAngleDependenceImpact}{\angleErosion} 
&:= (\sin{\angleErosion})^{\exponentBrittleOkaModel} (1+\hardnessOkaModel(1-\sin{\angleErosion}))^{\exponentDuctileOkaModel}, 
\\
 \dependence{\angleErosion}{\velocitySolid,\outerNormal} 
 &:= \frac{\pi}{2} - \arccos{\bigg(\frac{\transpose{\velocitySolid}\outerNormal}{\normEuclidean{\velocitySolid}} \bigg)}
  = \arcsin{\bigg(\frac{\transpose{\velocitySolid}\outerNormal}{\normEuclidean{\velocitySolid}} \bigg)},
\end{align*}
with $\exponentBrittleOkaModel,\exponentDuctileOkaModel,\hardnessOkaModel>0$
describe the impact angle-dependency.
The first factor of $\functionAngleDependenceImpact$ is associated with brittle material characteristics of the pipe, 
where the damaging is due to repeated plastic deformation at nearly orthogonal impact angles
and the second factor describes ductile material characteristics, 
where the tangential component of the particle velocity leads to small cuts at the surface of the pipe walls \cite{oka1997impact}.

Based on these definitions we consider the cost functional
\begin{equation}
 \label{eq:CostFunctionalEulerEulerModel}
  \dependence{\tilde{\costFunctionalLetter}}{\domain,\velocitySolid,\volumePercentageSolid}
  :=
 \IntegralSurfacePartDomainPlain{\dependence{\costFunctionalBoundaryFunction}{\volumePercentageSolid,\velocitySolid,\outerNormal}}{\surfaceCompDomainWall}
 + \dependence{\willmoreEnergy}{\domain}
\end{equation}
with $\dependence{\costFunctionalBoundaryFunction}{\volumePercentageSolid,\velocitySolid,\outerNormal}:= \frac{1}{2}\dependence{\functionalErosionLetter}{\volumePercentageSolid,\velocitySolid,\outerNormal}^2$,
where the squared erosion rate from \cref{eq:ErosionRateOKA} is used 
to penalize local maximums of the erosion rate at $\surfaceCompDomainWall$.
We note, that we do not take the material abrasion of the pipe over time 
for the geometrical description of $\domain$ into account,
so that the problem can be considered stationary.
The second integral in \cref{eq:CostFunctionalEulerEulerModel}, the so-called Willmore energy \cite{willmore1996riemannian} 
\begin{equation*}
 \dependence{\willmoreEnergy}{\domain}
 :=
 \constantRegularizationOne \IntegralSurfacePartDomainPlain{\frac{1}{2} \, \meanCurvature^2}{\surfaceCompDomainDeformable},
 \quad \constantRegularizationOne>0,
\end{equation*}
acts as a regularization penalizing surfaces with large mean curvature 
$\meanCurvature=\frac{1}{2} \tangentialDivergence{\outerNormal}$.

In the following we will assume, that \cref{eq:WeakForwardPDESystem}, 
the weak formulation of the Eulerian model,
possesses a unique solution,
so that the reduced cost functional 
\begin{equation}
\label{eq:CostFunctionReduced}
 \dependence{\costFunctionalLetter}{\domain} 
 := 
 \dependence{\tilde{\costFunctionalLetter}}{\domain,\dependence{\solutionSturmApproachVelocitySolid}{\domain},
 \dependence{\solutionSturmApproachVolumePercentageSolid}{\domain}}
  \end{equation}
can be defined, and consider the PDE-constrained optimization problem 
\begin{align}
\begin{split}
\label{eq:OptimizationProblem}
    \textup{Seek }
    \hat{\domain} \in \shapesAdmissable
    \textup{ such that } \\
      \dependence{\costFunctionalLetter}{\hat{\domain}} 
      = \min_{\domain \in \shapesAdmissable} 
	\dependence{\costFunctionalLetter}{\domain} 
	\\
    \textup{subject to } \eqref{eq:WeakForwardPDESystem} .
 \end{split}
 \end{align}

Here the set of admissible shapes $\shapesAdmissable \subset \powerset{\R^{3}}$ 
denotes the set of sufficiently smooth domains for which the non-deformable boundary part $\surfaceCompDomain \setminus \surfaceCompDomainDeformable$ is the same as for the initial domain, and we refer the reader to \cite[p.170]{DelfourZolesio2011shapes} for a more precise definition of $\shapesAdmissable$ in the general framework of shape optimization.

%% file: Pictures/tikzpicture_bend.tex
 \begin{tikzpicture}[scale=1.0]
  \tkzInit[xmax=4,ymax=4,xmin=-2,ymin=-3]
   \tkzDrawX[label=x]
   \tkzDrawY[label=z]
   \tkzLabelXY[text=black]
   
   \tikzstyle{ann} = [fill=white,font=\footnotesize,inner sep=0.5pt]
   
  \draw [red,domain=0:90] plot ({3.14*cos(\x)}, {3.14*sin(\x)});  
  \draw [red, domain=0:90] plot ({(3.14-2*0.395)*cos(\x)}, {(3.14-2*0.395)*sin(\x)});
    \node[ann] at (2.6,2.6) {\textcolor{red}{$\surfaceCompDomainDeformable \subset \surfaceCompDomainWall$}};  
  
  \draw [dash dot] (-1,3.14-0.395) -- (0,3.14-0.395); 
  \draw [domain=0:90,dash dot] plot ({(3.14-0.395)*cos(\x)}, {(3.14-0.395)*sin(\x)});
  \draw [dash dot] (3.14-0.395,0) -- (3.14-0.395,-2);   
  
  \draw [] (-1,3.14) -- (0,3.14); 
  \draw [] (3.14,0) -- (3.14,-2); 
  \draw [] (-1,3.14-2*0.395) -- (0,3.14-2*0.395); 
  \draw [] (3.14-2*0.395,0) -- (3.14-2*0.395,-2);
    \node[ann] at (3.14+0.3,-1) {$\surfaceCompDomainWall$};
    
  \draw [] (-1,3.14) -- (-1,3.14-2*0.395); 
    \node[ann] at (-1.3,3.14-0.395) {$\surfaceCompDomainInflow$};
  \draw [] (3.14-2*0.395,-2) -- (3.14,-2);
    \node[ann] at (3.14-0.395,-2.2) {$\surfaceCompDomainOutflow$};

  \draw [] (0,0) -- (1.9409,1.9409);
  \node[ann] at (1,1) {$\rCurvatureBend$};
  
  \draw[arrows=|-|](-1.8,3.14-2*0.395)--(-1.8,3.14);
  \node[ann] at (-1.8,3.14-0.395) {$\dTube$};

  \draw[arrows=|-|](-1,1.6)--(0,1.6);
  \node[ann] at (-0.5,1.6) {$\dTube$};

  \draw[arrows=|-|](1.5,0)--(1.5,-2);
    \node[ann] at (1.5,-1) {$2\dTube$};

  \draw[arrows=->,line width=0.25mm](4,2)--(4,1);
   \node[ann] at (3.8,1.5) {$\gravity$};
  
  \end{tikzpicture}

%% file: ShapeDerivatives.tex
The purpose of this section lies in the formal computation of the shape derivative of the reduced cost functional \cref{eq:CostFunctionReduced} in \cref{sec:ShapeDerivative}.
For examples of a more rigorous derivation of shape derivatives see e.g. \cite{DelfourZolesio2011shapes,gangl2015shape,hohmann2019residence,sturm2015minimax}.
In order to keep this chapter as self-contained as possible, 
we recall important definitions and results from shape calculus in \cref{sec:PreliminariesShapeCalculus}.
The continuous adjoint equations, 
which have to be solved in order to evaluate the shape derivative,
are derived in \cref{sec:AdjointEquations} through the introduction of the shape Lagrangian associated with the optimization problem \cref{eq:OptimizationProblem}.

%
%
\subsection{Preliminaries}
\label{sec:PreliminariesShapeCalculus}

We denote by 
\begin{equation*}
\functionSpacePertubationsAdmissableTwoRThree
:= \{\pertubationLetter \in \lipschitzContinuousFunctionsTwoRThree : \normLipschitzContinuousFunctionsTwoRThree{\pertubationLetter} < 0.5,
\,  \evaluateAtTwo{\pertubationLetter}{\surfaceCompDomain \setminus \surfaceCompDomainDeformable} = 0 
\}
\end{equation*}
the space of admissible deformations, equipped with the norm
\begin{equation}
\label{eq:NormOfC1Deformations}
 \normLipschitzContinuousFunctionsTwoRThree{\genericFunction}
 := \sup_{\abs{\multiindexDerivative}\le 2} 
 \normLInfty{\multiindexD{\genericFunction}{\multiindexDerivative}}
 + \sup_{\substack{\abs{\multiindexDerivative} = 2 \\
		  \spaceVar,\spaceVarTwo \in \domain \\
		  \spaceVar \neq \spaceVarTwo
		  }
	}\frac{\abs{\dependence{\multiindexD{\genericFunction}{\multiindexDerivative}}{\spaceVar} - \dependence{\multiindexD{\genericFunction}{\multiindexDerivative}}{\spaceVarTwo}}}{\abs{\spaceVar - \spaceVarTwo}},
	\quad \multiindexDerivative \in \N_{0}^{3}.
\end{equation}
Given $\domain \subset \R^{3}$
 and $\domainPertubationLetter \boundaryVelocity \in \functionSpacePertubationsAdmissableTwoRThree$ 
 for all $\domainPertubationLetter \in \closedInterval{0}{1}$,
the transformations
 \begin{align}
 \begin{split}
 \label{eq:PertubationOfTheIdentityApproach}
   \defineFunctionLong{\transformationFromReferenceDomain}{\domain}{\R^{3}}{\spaceVar}{\dependence{\transformationFromReferenceDomain}{\spaceVar}
   :=\dependence{(\identityFunction + \domainPertubationLetter\boundaryVelocity)}{\spaceVar}}
 \end{split}
\end{align}
 are diffeomorphisms due to the Neumann series \cite{simon1980differentiation} 
 and we set
$
  \domainControlled := \image{\transformationFromReferenceDomain}{\domain}.
 $
Note that the curvature can be defined on 
$\image{\transformationFromReferenceDomain}{\surfaceCompDomainDeformable}$, 
since it enjoys the same smoothness properties as $\surfaceCompDomainDeformable$ 
due to the assumed regularity of the admissible deformations.
 
 We further denote by $\jacTrafo = \identityMatrix + \domainPertubationLetter \jacDefVelo$ 
 the Jacobian of the transformation $\transformationFromReferenceDomain$ and by 
  \begin{subequations}
 \begin{align}
  \label{eq:InverseTransposedJacobianTrafo}
 \invTranspJacTrafoEvaluated  	&:= \inverseTranspose{\jacTrafo}, \\
  \label{eq:detTrafoDomain}
 \detTrafoEvaluated 		&:= \det \jacTrafo, \\
  \label{eq:detTrafoBoundary}
 \detTrafoBdryEvaluated  	&:= \detTrafoEvaluated \abs{\inverseTranspose{\jacTrafo}\outerNormal}, 
\end{align}
\end{subequations}
 the transposed inverse of this matrix, as well as the volume and surface elements of the transformation, which satisfy 
 $\dependence{\invTranspJacTrafo}{0} = \identityMatrix$
 and 
 $\dependence{\detTrafo}{0} = \dependence{\detTrafoBdry}{0} = 1$.
 The following Lemma,
 which will be applied in \Cref{thm:ShapeDerivative} for the computation of the shape derivative of 
 \cref{eq:CostFunctionReduced},
 addresses the derivatives of \cref{eq:InverseTransposedJacobianTrafo,eq:detTrafoDomain,eq:detTrafoBoundary},
 and a chain rule for the divergence and gradient composed with a smooth transformation of the domain given by \cref{eq:PertubationOfTheIdentityApproach}.
 \begin{lem}
  \label{lem:PropertiesOfTransformationAndDifferentialOperators}
  For $\genericFunction \in \spaceForGradientOfComposedScalarFunction$, 
  $\genericVectorFunction \in \spaceForGradientOfComposedVectorFunction$ and $\boundaryVelocity \in \lipschitzContinuousFunctionsOneRThree$ it holds, that 
\begin{equation}
\label{eq:DifferentiationOfTransformationQuantities}
 \dependence{\derivative{\detTrafo}}{0} 
 = \divergence{\boundaryVelocity} ,
 \quad 
 \dependence{\derivative{\detTrafoBdry}}{0} 
 = \tangentialDivergence{\boundaryVelocity},
 \quad 
 \dependence{\derivative{\invTranspJacTrafo}}{0} 
 = - \transpose{\jacDefVelo},
\end{equation}
with the tangential divergence 
$\tangentialDivergence{\boundaryVelocity}:= \divergence{\boundaryVelocity} - \transpose{\outerNormal}\jacobian{\boundaryVelocity} \outerNormal$.

The outer normal $\outerNormalTransformed$ on $\image{\transformationFromReferenceDomain}{\surfaceCompDomain}$
satisfies 
\begin{align}
  \label{eq:TotalDerivativeOuterNormalOnTransformedDomainWithPushForward}
  \evaluateAtTwo{ \totalDNoFrac{(\outerNormalTransformed \circ \transformationFromReferenceDomain)}{\domainPertubationLetter} }{\domainPertubationLetter=0}
  = \evaluateAtTwo{ \totalDNoFrac{
			   \bigg( 
			   \frac{\invTranspJacTrafoEvaluated \, \outerNormal}{\normEuclidean{\invTranspJacTrafoEvaluated \, \outerNormal}}
			  \bigg)
		    }{\domainPertubationLetter} }{\domainPertubationLetter=0}
  = - \transpose{(\tangentialJacobian{\boundaryVelocity})} \outerNormal
\end{align}
with the tangential jacobian 
$ \tangentialJacobian{\boundaryVelocity} 
    := \jacobian{\boundaryVelocity} - \jacobian{\boundaryVelocity} \outerNormal \transpose{\outerNormal}$.
Further, 
    \begin{subequations}
    \begin{align}
    \label{eq:TransformationOfDifferentialOperatorsOne}
    \gradient{(\genericFunction \circ \inverse{\transformationFromReferenceDomain})}
    &= (\invTranspJacTrafoEvaluated \, \gradient{\genericFunction} )
    \circ \inverse{\transformationFromReferenceDomain}, \\
    \label{eq:TransformationOfDifferentialOperatorsTwo}
      \divergence{(\genericVectorFunction \circ \inverse{\transformationFromReferenceDomain})}
    &= \trace{\jacobian{\genericVectorFunction} \, \invJacTrafoEvaluated}\circ \inverse{\transformationFromReferenceDomain}.
    \end{align}
    \end{subequations}
\end{lem}

 \begin{pf}
  The proofs of \cref{eq:DifferentiationOfTransformationQuantities} and the first equality in
  \cref{eq:TransformationOfDifferentialOperatorsOne,eq:TransformationOfDifferentialOperatorsTwo} can be found in \cite{DelfourZolesio2011shapes}
  and \cref{eq:TotalDerivativeOuterNormalOnTransformedDomainWithPushForward}
  is shown in \cite{Schmidt2010shape}.
  Furthermore
  \begin{equation*}
  \divergence{(\genericVectorFunction \circ \inverse{\transformationFromReferenceDomain})}
  =
 \trace{\jacobian{(\genericVectorFunction\circ \inverse{\transformationFromReferenceDomain})}}
 = \trace{(\jacobian{\genericVectorFunction}\circ \inverse{\transformationFromReferenceDomain}) \, \inverse{\jacobian{\transformationFromReferenceDomain}}}
 = \trace{(\jacobian{\genericVectorFunction} \, \invJacTrafoEvaluated)} \circ \inverse{\transformationFromReferenceDomain}
\end{equation*}
holds due to the chain rule.
 \end{pf}

 We now define derivatives of shape functionals with respect to changes of the domain.
 \begin{defn}[\cite{Sokolowski1992introduction}]
  \label{def:EulerianDerivativeOfCostFunctional}
 $\defineFunctionShortSpacesOnly{\costFunctionalLetter}{\shapesAdmissable}{\R}$ is said to have the Eulerian semi-derivative in direction $\boundaryVelocity \in \lipschitzContinuousFunctionsGeneralRThree$ for $\orderShapeDerivative \in \N$,  if the limit
 \begin{equation}
  \label{eq:DefinitionShapeDerivative}
  \denoteEulerianDerivativeFunctionalDomain{\costFunctionalLetter}{\boundaryVelocity}
  := \lim_{\domainPertubationLetter \searrow 0} 
  \frac{ 
	      \dependence{\costFunctionalLetter}{\dependence{\transformationFromReferenceDomain}{\domain}} - \dependence{\costFunctionalLetter}{\domain}
	  }
	{\domainPertubationLetter}  
 \end{equation}
exists.
It is called shape differentiable of order $\orderShapeDerivative$, if this limit exists for all $\boundaryVelocity \in \lipschitzContinuousFunctionsGeneralRThree$ and the mapping 
$\boundaryVelocity \mapsto \denoteEulerianDerivativeFunctionalDomain{\costFunctionalLetter}{\boundaryVelocity}$
 is linear and continuous.
 \Cref{eq:DefinitionShapeDerivative} is then called shape-derivative of $J$.
 \end{defn}

\subsection{Adjoint equations}
\label{sec:AdjointEquations}
 
%
%

In order to compute the Eulerian semi-derivative of \cref{eq:CostFunctionReduced},
we have to determine the corresponding adjoint equations.
This is done by introducing the shape Lagrangian 
associated with the optimization problem
and following the approach described in \cite[Chap. 10]{DelfourZolesio2011shapes}.
Given 
$\ansatzSpaceSturmApproachLetter
:=(\ansatzSpaceSturmApproachVelocity,\ansatzSpaceSturmApproachPressure,\ansatzSpaceSturmApproachVelocitySolid, \ansatzSpaceSturmApproachVolumePercentageSolid)
\in \spaceE$
and 
$\testSpaceSturmApproachLetter:=
(\testSpaceSturmApproachVelocity, \testSpaceSturmApproachPressure, \testSpaceSturmApproachVelocitySolid, \testSpaceSturmApproachVolumePercentageSolid, \testSpaceSturmApproachVelocityFluidInflowCondition, \testSpaceSturmApproachVelocitySolidInflowCondition, \testSpaceSturmApproachVolumePercentageSolidInflowCondition) 
\in \spaceF$ 
we introduce the PDE-constraint 
\begin{align*}
 \dependence{\functionalEq}{\domain,\ansatzSpaceSturmApproachLetter,\testSpaceSturmApproachLetter} 
 :=&
  \dependence{\functionalC}{\domain,\ansatzSpaceSturmApproachVelocity,\ansatzSpaceSturmApproachVelocity,\testSpaceSturmApproachVelocity}
  + \reciprocal{\reynoldsNumber}\dependence{\functionalA}{\domain,\ansatzSpaceSturmApproachVelocity,\testSpaceSturmApproachVelocity} 
  - \dependence{\functionalB}{\domain,\testSpaceSturmApproachVelocity,\ansatzSpaceSturmApproachPressure}
 - {\froudeNumber}^{-2}\dependence{\functionalG}{\domain,\testSpaceSturmApproachVelocity} \\
 &+ \dependence{\functionalB}{\domain,\ansatzSpaceSturmApproachVelocity,\testSpaceSturmApproachPressure}
 + \dependence{\functionalC}{\domain,\ansatzSpaceSturmApproachVelocitySolid,\ansatzSpaceSturmApproachVelocitySolid,\testSpaceSturmApproachVelocitySolid}
 +\reciprocal{\smoothingVelocitySolidNumber}
   \dependence{\functionalA}{\domain,\ansatzSpaceSturmApproachVelocitySolid,\testSpaceSturmApproachVelocitySolid} \\
   &+ 2\,\reciprocal{\stokesNumber} \dependence{\functionalF}{\domain, \ansatzSpaceSturmApproachVelocity,\ansatzSpaceSturmApproachVelocitySolid,\testSpaceSturmApproachVelocitySolid}
   - {\froudeNumber}^{-2}\dependence{\functionalG}{\domain,\testSpaceSturmApproachVelocitySolid} \\
   &+  \dependence{\functionalD}{\domain,\ansatzSpaceSturmApproachVolumePercentageSolid,\ansatzSpaceSturmApproachVelocitySolid,\testSpaceSturmApproachVolumePercentageSolid}
   +\reciprocal{\pecletNumber}\dependence{\functionalA}{\ansatzSpaceSturmApproachVolumePercentageSolid,\testSpaceSturmApproachVolumePercentageSolid} \\
   &+\dependence{\functionalHVelocityFluidInflow}{\domain,\ansatzSpaceSturmApproachVelocity,\testSpaceSturmApproachVelocityFluidInflowCondition}
   +\dependence{\functionalHVelocitySolidInflow}{\domain,\ansatzSpaceSturmApproachVelocitySolid,\testSpaceSturmApproachVelocitySolidInflowCondition}
   +\dependence{\functionalHVolumePercentageSolidInflow}{\domain,\ansatzSpaceSturmApproachVolumePercentageSolid,\testSpaceSturmApproachVolumePercentageSolidInflowCondition}
\end{align*}
of the weak formulation \cref{eq:WeakForwardPDESystem}.
The shape Lagrangian associated with \cref{eq:OptimizationProblem} is defined through
\begin{equation}
\label{eq:ShapeLagrangianOnDomain}
 \dependence{\lagrangianSturmApproach}{\domain,\ansatzSpaceSturmApproachLetter,\testSpaceSturmApproachLetter}
 := \dependence{\tilde{\costFunctionalLetter}}{\domain,\ansatzSpaceSturmApproachVelocitySolid,\ansatzSpaceSturmApproachVolumePercentageSolid}
  +  \dependence{\functionalEq}{\domain,\ansatzSpaceSturmApproachLetter,\testSpaceSturmApproachLetter},
\end{equation}
and we note that 
the weak solution $\solutionSturmApproachDomain$ of \cref{eq:WeakForwardPDESystem}
satisfies 
\begin{equation*}
 \directionalDHilbertSpaces{\testSpaceSturmApproachLetter}{\dependence{\lagrangianSturmApproach}{\domain,\solutionSturmApproachDomain,\cdot}}\testSpaceSturmApproachLetter
 = 0 \quad \forall \, \testSpaceSturmApproachLetter \in \spaceF.
\end{equation*}

The adjoint state
$\adjointSolutionSturmApproachDomain:= (\adjointSolutionSturmApproachVelocityFluid,\adjointSolutionSturmApproachPressure,\adjointSolutionSturmApproachVelocitySolid,\adjointSolutionSturmApproachVolumePercentageSolid,\adjointSolutionSturmApproachVelocityFluidInflowCondition,\adjointSolutionSturmApproachVelocitySolidInflowCondition,\adjointSolutionSturmApproachVolumePercentageSolidInflowCondition) \in\spaceF$
is defined as the solution of
\begin{equation}
\label{eq:AdjointEquationFromShapeLagrangian}
 \directionalDHilbertSpaces{\ansatzSpaceSturmApproachLetter}{\dependence{\lagrangianSturmApproach}{\domain,\solutionSturmApproachDomain,\adjointSolutionSturmApproachDomain}}\ansatzSpaceSturmApproachLetter
 =0 \quad \forall \, \ansatzSpaceSturmApproachLetter \in \spaceE,
\end{equation}
and in order to explicitly derive the adjoint equations, 
we take the derivative with respect to the individual components of $\ansatzSpaceSturmApproachLetter$.

Choosing $\ansatzSpaceSturmApproachLetter=(0,0,0,\ansatzSpaceSturmApproachVolumePercentageSolid)$ in \cref{eq:AdjointEquationFromShapeLagrangian}
and using integration by parts, 
we obtain 
the 
adjoint volume percentage transport equation 
\begin{align}
 \begin{split}
 \label{eq:AdjointEquationParticleVolumePercentageEulerEulerModel}
 &\IntegralSurfacePartDomainPlain{
 \ansatzSpaceSturmApproachVolumePercentageSolid \adjointSolutionSturmApproachVolumePercentageSolidInflowCondition }{\surfaceCompDomainInflow}
 +\IntegralDomain{ -\transpose{\solutionSturmApproachVelocitySolid} \, \gradient{\adjointSolutionSturmApproachVolumePercentageSolid} \,\ansatzSpaceSturmApproachVolumePercentageSolid }
  + \reciprocal{\pecletNumber}
  \IntegralDomain{\transpose{\gradient{\adjointSolutionSturmApproachVolumePercentageSolid}}\,
  \gradient{\ansatzSpaceSturmApproachVolumePercentageSolid}} 
 \\
  &+ \IntegralSurfacePartDomainPlain{
  \ansatzSpaceSturmApproachVolumePercentageSolid (\transpose{\solutionSturmApproachVelocitySolid}  \outerNormal)
  \adjointSolutionSturmApproachVolumePercentageSolid 
  }{\dependence{\surfaceCompDomainEndingCharacteristics}{\solutionSturmApproachVelocitySolid}}
 =
  -  \IntegralSurfacePartDomainPlain{\dependence{\partialD{\costFunctionalBoundaryFunction}{\ansatzSpaceSturmApproachVolumePercentageSolid}}{\solutionSturmApproachVolumePercentageSolid,\solutionSturmApproachVelocitySolid,\outerNormal}\ansatzSpaceSturmApproachVolumePercentageSolid
  }{\surfaceCompDomainWall}
  \quad \forall \,  \ansatzSpaceSturmApproachVolumePercentageSolid \in \spaceZ .
\end{split}
\end{align}

Choosing $\ansatzSpaceSturmApproachLetter=(0,0,\ansatzSpaceSturmApproachVelocitySolid,0)$
in \cref{eq:AdjointEquationFromShapeLagrangian},
we derive
\begin{align}
\begin{split}
\label{eq:AdjointEquationParticleVelocityEulerEulerModel}
   &  \IntegralSurfacePartDomainPlain{\transpose{\adjointSolutionSturmApproachVelocitySolidInflowCondition} \ansatzSpaceSturmApproachVelocitySolid }{\surfaceCompDomainInflow}
   +\IntegralDomain{
   (\transpose{\adjointSolutionSturmApproachVelocitySolid} \, \gradient{\solutionSturmApproachVelocitySolid} 
  - \transpose{\solutionSturmApproachVelocitySolid}\, \transpose{\gradient{\adjointSolutionSturmApproachVelocitySolid}}
  - (\divergence{\solutionSturmApproachVelocitySolid}) \transpose{\adjointSolutionSturmApproachVelocitySolid} \, )\ansatzSpaceSturmApproachVelocitySolid }  
   + \IntegralDomain{\divergence{(\solutionSturmApproachVolumePercentageSolid \ansatzSpaceSturmApproachVelocitySolid)}\adjointSolutionSturmApproachVolumePercentageSolid} 
 \\
  &+ 2\reciprocal{\stokesNumber} 
      \IntegralDomain{\transpose{(\dependence{\partialD{\weakContributionSchillerNaumannDragTerm}{\ansatzSpaceSturmApproachVelocitySolid}}{\solutionSturmApproachVelocityFluid,\solutionSturmApproachVelocitySolid,\adjointSolutionSturmApproachVelocitySolid})} \ansatzSpaceSturmApproachVelocitySolid} 
  + \reciprocal{\smoothingVelocitySolidNumber} \IntegralDomain{  \gradient{\adjointSolutionSturmApproachVelocitySolid}:
  \gradient{\ansatzSpaceSturmApproachVelocitySolid}} 
    + \IntegralSurfacePartDomainPlain{(\transpose{\solutionSturmApproachVelocitySolid} \outerNormal) \transpose{\adjointSolutionSturmApproachVelocitySolid}  \ansatzSpaceSturmApproachVelocitySolid}{\surfaceCompDomain \setminus \surfaceCompDomainInflow} \\
  &=  
  - \IntegralSurfacePartDomainPlain{\transpose{\dependence{\partialD{\costFunctionalBoundaryFunction}{\ansatzSpaceSturmApproachVelocitySolid}}{\solutionSturmApproachVolumePercentageSolid,\solutionSturmApproachVelocitySolid,\outerNormal}}\ansatzSpaceSturmApproachVelocitySolid
  }{\surfaceCompDomainWall}
  \quad \forall \, \ansatzSpaceSturmApproachVelocitySolid \in \spaceY,
\end{split}
\end{align}
where the contribution 
\begin{equation*}
\dependence{\weakContributionSchillerNaumannDragTerm}{\ansatzSpaceSturmApproachVelocity,\ansatzSpaceSturmApproachVelocitySolid}
:= \dependence{\SchillerNaumannDragTerm}{\ansatzSpaceSturmApproachVelocity,\ansatzSpaceSturmApproachVelocitySolid} \transpose{(\ansatzSpaceSturmApproachVelocitySolid - \ansatzSpaceSturmApproachVelocity)} \adjointSolutionSturmApproachVelocitySolid
\end{equation*}
stems from the drag term. Its derivatives are given through
\begin{align*}
      \transpose{(
      \dependence{\partialD{\weakContributionSchillerNaumannDragTerm}{\ansatzSpaceSturmApproachVelocity}}{\solutionSturmApproachVelocityFluid,\solutionSturmApproachVelocitySolid}
      )
      }
      \ansatzSpaceSturmApproachVelocity
      &= (\transpose{(\solutionSturmApproachVelocitySolid-\solutionSturmApproachVelocityFluid)}\adjointSolutionSturmApproachVelocitySolid)
      \transpose{
      \dependence{\partialD{\SchillerNaumannDragTerm}{\ansatzSpaceSturmApproachVelocity}}{\solutionSturmApproachVelocityFluid,\solutionSturmApproachVelocitySolid}
      }\ansatzSpaceSturmApproachVelocity
      - \dependence{\SchillerNaumannDragTerm}{\solutionSturmApproachVelocityFluid,\solutionSturmApproachVelocitySolid}
	    \transpose{\adjointSolutionSturmApproachVelocitySolid} \ansatzSpaceSturmApproachVelocity, \\
\transpose{(
\dependence{\partialD{\weakContributionSchillerNaumannDragTerm}{\ansatzSpaceSturmApproachVelocitySolid}}{\solutionSturmApproachVelocityFluid,\solutionSturmApproachVelocitySolid}
 )}
 \ansatzSpaceSturmApproachVelocitySolid
 &= (\transpose{(\solutionSturmApproachVelocitySolid-\solutionSturmApproachVelocityFluid)}\adjointSolutionSturmApproachVelocitySolid)
 \transpose{
 \dependence{\partialD{\SchillerNaumannDragTerm}{\ansatzSpaceSturmApproachVelocitySolid}}{\solutionSturmApproachVelocityFluid,\solutionSturmApproachVelocitySolid}
 }\ansatzSpaceSturmApproachVelocitySolid
 + \dependence{\SchillerNaumannDragTerm}{\solutionSturmApproachVelocityFluid,\solutionSturmApproachVelocitySolid}
      \transpose{\adjointSolutionSturmApproachVelocitySolid} \ansatzSpaceSturmApproachVelocitySolid \\
&=-\transpose{(
      \dependence{\partialD{\weakContributionSchillerNaumannDragTerm}{\ansatzSpaceSturmApproachVelocity}}{\solutionSturmApproachVelocityFluid,\solutionSturmApproachVelocitySolid}
      )
      }
      \ansatzSpaceSturmApproachVelocitySolid,
\end{align*}
with 
\begin{equation*}
\dependence{\partialD{\SchillerNaumannDragTerm}{\ansatzSpaceSturmApproachVelocityComponent}}{\solutionSturmApproachVelocityFluid,\solutionSturmApproachVelocitySolid}
=  
 0.10305 \, (\densityFluid \diamEq \reciprocal{\viscosityFluid})^{0.687}
  \frac{\solutionSturmApproachVelocityFluidComponent-\solutionSturmApproachVelocitySolidComponent}{\normEuclidean{\solutionSturmApproachVelocitySolid-\solutionSturmApproachVelocityFluid}^{1.313}},
  \qquad 
    \dependence{\partialD{\SchillerNaumannDragTerm}{\ansatzSpaceSturmApproachVelocitySolidComponent}}{\solutionSturmApproachVelocityFluid,\solutionSturmApproachVelocitySolid}
  = -
    \dependence{\partialD{\SchillerNaumannDragTerm}{\ansatzSpaceSturmApproachVelocityComponent}}{\solutionSturmApproachVelocityFluid,\solutionSturmApproachVelocitySolid}.
\end{equation*}

The adjoint fluid velocity equation
\begin{align}
\begin{split}
\label{eq:AdjointEquationFluidVelocityEulerEulerModel}
&
 \IntegralSurfacePartDomainPlain{\transpose{\adjointSolutionSturmApproachVelocityFluidInflowCondition} \ansatzSpaceSturmApproachVelocity }{\surfaceCompDomainInflow}
 +\IntegralDomain{ \transpose{\adjointSolutionSturmApproachVelocityFluid} \, \gradient{\solutionSturmApproachVelocityFluid} \,
\ansatzSpaceSturmApproachVelocity
 - \transpose{\solutionSturmApproachVelocityFluid}\, \transpose{\gradient{\adjointSolutionSturmApproachVelocityFluid}}\, \ansatzSpaceSturmApproachVelocity}
  + \IntegralSurfacePartDomainPlain{(\transpose{\solutionSturmApproachVelocityFluid} \outerNormal) (\transpose{\adjointSolutionSturmApproachVelocityFluid }  \ansatzSpaceSturmApproachVelocity)}{\surfaceCompDomainOutflow} \\
  &+ \reciprocal{\reynoldsNumber} \IntegralDomain{ \gradient{\adjointSolutionSturmApproachVelocityFluid}:
  \gradient{\ansatzSpaceSturmApproachVelocity}} 
  + \IntegralDomain{\adjointSolutionSturmApproachPressure \divergence{\ansatzSpaceSturmApproachVelocity}} 
  \\
  &=  
  2\reciprocal{\stokesNumber}
      \IntegralDomain{\transpose{(\dependence{\partialD{\weakContributionSchillerNaumannDragTerm}{\ansatzSpaceSturmApproachVelocity}}{\solutionSturmApproachVelocityFluid,\solutionSturmApproachVelocitySolid,\adjointSolutionSturmApproachVelocitySolid})} \, \ansatzSpaceSturmApproachVelocity} 
  \quad \forall \, \ansatzSpaceSturmApproachVelocity \in \spaceX,
\end{split}
\end{align}
and the incompressibility condition
\begin{equation}
\label{eq:AdjointEquationIncompressibilityEulerEulerModel} 
\IntegralDomain{\ansatzSpaceSturmApproachPressure \divergence{\adjointSolutionSturmApproachVelocityFluid}}
 = 0 \quad \forall \, \ansatzSpaceSturmApproachPressure \in \spaceM,
\end{equation}
are obtained by choosing 
$\ansatzSpaceSturmApproachLetter=(\ansatzSpaceSturmApproachVelocity,0,0,0)$
and 
$\ansatzSpaceSturmApproachLetter=(0,\ansatzSpaceSturmApproachPressure,0,0)$
respectively in \cref{eq:AdjointEquationFromShapeLagrangian}.

We note, that due to the decoupling of the forward equations, 
the components of the adjoint state $\adjointSolutionSturmApproachDomain$ can also be obtained 
in a decoupled manner 
by solving the sub-problems
\cref{eq:AdjointEquationParticleVolumePercentageEulerEulerModel}, \cref{eq:AdjointEquationParticleVelocityEulerEulerModel} as well as 
\cref{eq:AdjointEquationFluidVelocityEulerEulerModel,eq:AdjointEquationIncompressibilityEulerEulerModel}
in order.
Even though we restrict ourselves to the OKA erosion model \cite{oka2005practical}, 
the adjoint transport equation \cref{eq:AdjointEquationParticleVolumePercentageEulerEulerModel}
and the adjoint particle velocity equation \cref{eq:AdjointEquationParticleVelocityEulerEulerModel}
remain valid for the volume-averaged formulation of the erosion model of Finnie and the E/CRC model,
since only the partial derivatives of $\costFunctionalBoundaryFunction$ have to be replaced accordingly.

\subsection{Shape derivative}
\label{sec:ShapeDerivative}

Given the solutions of the forward and adjoint problem on $\domain$, 
we formally compute the shape derivative of \cref{eq:CostFunctionReduced} in the following
theorem.
\begin{thm}
\label{thm:ShapeDerivative}
Let $\boundaryVelocity \in \lipschitzContinuousFunctionsTwoRThree$. 
If it exists, 
the Eulerian semi-derivative of \cref{eq:CostFunctionReduced} in direction $\boundaryVelocity$ is given through
  \begin{align}
      \label{eq:ShapeDerivativeEulerEulerModel}
  \denoteEulerianDerivativeFunctionalDomain{\costFunctionalLetter}{\boundaryVelocity} 
     &=  \IntegralSurfacePartDomainPlain{
	      \dependence{\costFunctionalBoundaryFunction}{\solutionSturmApproachVolumePercentageSolid,\solutionSturmApproachVelocitySolid,\outerNormal} \tangentialDivergence{\boundaryVelocity}
	      -\transpose{(\dependence{\partialD{\costFunctionalBoundaryFunction}{\outerNormal}}{\solutionSturmApproachVolumePercentageSolid,\solutionSturmApproachVelocitySolid,\outerNormal})}\,
	      \transpose{(\tangentialJacobian{\boundaryVelocity})} \outerNormal
	}{\surfaceCompDomainWall}  \nonumber \\
  &+ \constantRegularizationOne
      \IntegralSurfacePartDomainPlain{ 
	\big(
	  (\identityMatrix - (\tangentialJacobian{\spaceVar}+\transpose{(\tangentialJacobian{\spaceVar})}))
	  \,
	  \tangentialJacobian{\boundaryVelocity} 
	 \big)
	     :
	     \big(\tangentialJacobian{(\meanCurvature \outerNormal)}\big)
	+ \frac{1}{2}\, (\tangentialDivergence{(\meanCurvature \outerNormal)}) (\tangentialDivergence{\boundaryVelocity})      	      
      }{\surfaceCompDomainDeformable} \nonumber \\
 &+ \IntegralDomain{
      \big(
      \transpose{
      (\solutionSturmApproachVelocityFluid} \transpose{\gradient{\solutionSturmApproachVelocityFluid}}
       - \froudeNumber^{-2} \transpose{\gravity} 
       )
       \adjointSolutionSturmApproachVelocityFluid
      +\reciprocal{\reynoldsNumber}\gradient{\solutionSturmApproachVelocityFluid}:\gradient{\adjointSolutionSturmApproachVelocityFluid}
      - \solutionSturmApproachPressure \divergence{\adjointSolutionSturmApproachVelocityFluid}
      + \adjointSolutionSturmApproachPressure \divergence{\solutionSturmApproachVelocityFluid}
      \big) \divergence{\boundaryVelocity}
    }   \nonumber \\
 &+ \IntegralDomain{
      \Big(
	\big(
	\transpose{\solutionSturmApproachVelocitySolid}  \transpose{\gradient{\solutionSturmApproachVelocitySolid}} 
	+2 \reciprocal{\stokesNumber} 
	\dependence{\SchillerNaumannDragTerm}{\solutionSturmApproachVelocityFluid,\solutionSturmApproachVelocitySolid}
	\transpose{(\solutionSturmApproachVelocitySolid - \solutionSturmApproachVelocityFluid)} 	
	- {\froudeNumber}^{-2} \transpose{\gravity}
	\big)
	\adjointSolutionSturmApproachVelocitySolid
      \Big)
      \divergence{\boundaryVelocity} 
 }  \nonumber \\
  &+ \IntegralDomain{
	\reciprocal{\smoothingVelocitySolidNumber} 
	(\gradient{\solutionSturmApproachVelocitySolid} : \gradient{\adjointSolutionSturmApproachVelocitySolid}
	)
	\divergence{\boundaryVelocity} 
	}
  +\IntegralDomain{
      \big(
	\adjointSolutionSturmApproachVolumePercentageSolid 
	\divergence{(\solutionSturmApproachVolumePercentageSolid \solutionSturmApproachVelocitySolid)}
	+\reciprocal{\pecletNumber}
	\transpose{\gradient{\solutionSturmApproachVolumePercentageSolid}}  \gradient{\adjointSolutionSturmApproachVolumePercentageSolid}
      \big)
      \divergence{\boundaryVelocity} 
 }  \\
 &- \IntegralDomain{
      \transpose{\solutionSturmApproachVelocityFluid}  \transpose{\jacobian{\boundaryVelocity}} \transpose{\gradient{\solutionSturmApproachVelocityFluid}}  \adjointSolutionSturmApproachVelocityFluid
      +\reciprocal{\reynoldsNumber}      
	  \big(
	    (\transpose{\jacobian{\boundaryVelocity}} \transpose{\gradient{\solutionSturmApproachVelocityFluid}})):\gradient{\adjointSolutionSturmApproachVelocityFluid}
	    +\gradient{\solutionSturmApproachVelocityFluid}
	      : (\transpose{\jacobian{\boundaryVelocity}}
	      \transpose{\gradient{\adjointSolutionSturmApproachVelocityFluid}})
	    \big)  
  }     \nonumber  \\
  & - \IntegralDomain{
	  \adjointSolutionSturmApproachPressure \trace{\jacobian{\solutionSturmApproachVelocityFluid}\jacobian{\boundaryVelocity}}
	  -\solutionSturmApproachPressure \trace{\jacobian{\adjointSolutionSturmApproachVelocityFluid}\jacobian{\boundaryVelocity}}
      } \nonumber \\
   & - \IntegralDomain{
	  \transpose{\solutionSturmApproachVelocitySolid} \transpose{\jacobian{\boundaryVelocity}} \transpose{\gradient{\solutionSturmApproachVelocitySolid}} \adjointSolutionSturmApproachVelocitySolid
	  +\reciprocal{\smoothingVelocitySolidNumber}
	  \big(
	    (\transpose{\jacobian{\boundaryVelocity}} \transpose{\gradient{\solutionSturmApproachVelocitySolid}})):\gradient{\adjointSolutionSturmApproachVelocitySolid}
	    +\gradient{\solutionSturmApproachVelocitySolid}
	      : (\transpose{\jacobian{\boundaryVelocity}}
	      \transpose{\gradient{\adjointSolutionSturmApproachVelocitySolid}})
	    \big)
	} \nonumber  \\
  & - \IntegralDomain{
	  \adjointSolutionSturmApproachVolumePercentageSolid 
	  \trace{\jacobian{(\solutionSturmApproachVolumePercentageSolid \solutionSturmApproachVelocitySolid)} \jacobian{\boundaryVelocity}}
    + \reciprocal{\pecletNumber}
	  \big(
	    (\transpose{\jacobian{\boundaryVelocity}} \gradient{\solutionSturmApproachVolumePercentageSolid}
	    ):
	    \gradient{\adjointSolutionSturmApproachVolumePercentageSolid}
	    +\gradient{\solutionSturmApproachVolumePercentageSolid} :(\transpose{\jacobian{\boundaryVelocity}}\gradient{\adjointSolutionSturmApproachVolumePercentageSolid})
	  \big)  
	}. \nonumber 
\end{align}
\end{thm}
\begin{pf}
Let 
\begin{equation}
  \label{eq:DefinitionPulledBackShapeLagrangian}
   \dependence{\lagrangianSturmApproachTransformed}{\domainPertubationLetter,\ansatzSpaceSturmApproachLetter,\testSpaceSturmApproachLetter} 
   := \dependence{\lagrangianSturmApproach}{\domainControlled,\ansatzSpaceSturmApproachLetter \circ \inverse{\transformationFromReferenceDomain},\testSpaceSturmApproachLetter \circ \inverse{\transformationFromReferenceDomain}}
 \end{equation}
with the shape Lagrangian from \cref{eq:ShapeLagrangianOnDomain}. 
In \cite{DelfourZolesio2011shapes} it is shown, 
that if the Eulerian semi-derivative 
of \cref{eq:CostFunctionReduced} exists, 
it can be computed through 
\begin{equation}
  \label{eq:ComputationOfEulerianSemiderivative}
 \denoteEulerianDerivativeFunctionalDomain{\costFunctionalLetter}{\boundaryVelocity}
 = \dependence{\partialD{\lagrangianSturmApproachTransformed}{\domainPertubationLetter}}{0,\solutionSturmApproachDomain,\adjointSolutionSturmApproachDomain},
\end{equation}
where $\solutionSturmApproachDomain\in\spaceE$ 
if the solution of the weak forward problem \cref{eq:WeakForwardPDESystem} and 
$\adjointSolutionSturmApproachDomain\in\spaceF$ 
the solution of the weak adjoint problem \cref{eq:AdjointEquationFromShapeLagrangian}.
The domain-dependency of the function spaces 
of the ansatz and test functions 
$\spaceEDeformed$ and $\spaceFDeformed$ in \cref{eq:DefinitionPulledBackShapeLagrangian}
 is circumvented by parametrizing
 $\spaceEDeformed$ and $\spaceFDeformed$ by elements of $\spaceE$ and $\spaceF$ composed with $\inverse{\transformationFromReferenceDomain}$, see e.g.
 \cite[Theorem 2.2.2, p. 52]{ziemer2012weakly}.
In order to compute the partial derivative on the right hand side of \cref{eq:ComputationOfEulerianSemiderivative},
we use the transformation rule to pull back the integrals in \cref{eq:DefinitionPulledBackShapeLagrangian} to $\domain$,
as well as  \cref{eq:TransformationOfDifferentialOperatorsOne,eq:TransformationOfDifferentialOperatorsTwo}.
Since the derivative of the Willmore functional is given in \cite{bonito2010parametric}, we consider
\begingroup
  \allowdisplaybreaks
\begin{align}
   \label{eq:TransformedLagrangianEulerEulerModel}
 &\dependence{\lagrangianSturmApproachTransformed}{\domainPertubationLetter,\ansatzSpaceSturmApproachLetter,\testSpaceSturmApproachLetter} 
 -  \dependence{\willmoreEnergy}{\domainControlled}
  \nonumber \\
 &=
    \IntegralSurfacePartDomainPlain{\dependence{\costFunctionalBoundaryFunction}{\ansatzSpaceSturmApproachVolumePercentageSolid,\ansatzSpaceSturmApproachVelocitySolid,\outerNormalTransformed \circ \transformationFromReferenceDomain} \detTrafoBdryEvaluated}{\surfaceCompDomainWall} 
      +\IntegralDomain{\transpose{\ansatzSpaceSturmApproachVelocity} \, \invTranspJacTrafoEvaluated  \, \transpose{\gradient{\ansatzSpaceSturmApproachVelocity}} \, \testSpaceSturmApproachVelocity \detTrafoEvaluated}  \nonumber \\
      &+ \reciprocal{\reynoldsNumber}\IntegralDomain{(\invTranspJacTrafoEvaluated  \transpose{\gradient{\ansatzSpaceSturmApproachVelocity}}):(\invTranspJacTrafoEvaluated  \transpose{\gradient{\testSpaceSturmApproachVelocity}}) \detTrafoEvaluated}
      - {\froudeNumber}^{-2} \IntegralDomain{\transpose{\gravity}  \testSpaceSturmApproachVelocity \detTrafoEvaluated} \nonumber \\
      & -\IntegralDomain{ \ansatzSpaceSturmApproachPressure \trace{\jacobian{\testSpaceSturmApproachVelocity} \, \invJacTrafoEvaluated}\detTrafoEvaluated} 
      +\IntegralDomain{\testSpaceSturmApproachPressure \trace{\jacobian{\ansatzSpaceSturmApproachVelocity} \,\invJacTrafoEvaluated}\detTrafoEvaluated} \nonumber \\
      &+\IntegralDomain{\transpose{\ansatzSpaceSturmApproachVelocitySolid}\, \invTranspJacTrafoEvaluated  \, \transpose{\gradient{\ansatzSpaceSturmApproachVelocitySolid}}\, \testSpaceSturmApproachVelocitySolid \detTrafoEvaluated} 
    + 2\,\reciprocal{\stokesNumber} \IntegralDomain{
    \dependence{\SchillerNaumannDragTerm}{\ansatzSpaceSturmApproachVelocity,\ansatzSpaceSturmApproachVelocitySolid} \,
    \transpose{(\ansatzSpaceSturmApproachVelocitySolid-\ansatzSpaceSturmApproachVelocity)} \,  \testSpaceSturmApproachVelocitySolid \detTrafoEvaluated}  \\
    &+  \reciprocal{\smoothingVelocitySolidNumber}\IntegralDomain{(\invTranspJacTrafoEvaluated  \transpose{\gradient{\ansatzSpaceSturmApproachVelocitySolid}}):(\invTranspJacTrafoEvaluated \transpose{\gradient{\testSpaceSturmApproachVelocitySolid}}) \detTrafoEvaluated}
    - {\froudeNumber}^{-2} \IntegralDomain{\transpose{\gravity}  \testSpaceSturmApproachVelocitySolid \detTrafoEvaluated} \nonumber \\ 
    &+ 
    \IntegralDomain{
    \testSpaceSturmApproachVolumePercentageSolid
    \trace{\jacobian{(\ansatzSpaceSturmApproachVolumePercentageSolid \ansatzSpaceSturmApproachVelocitySolid)}\invJacTrafoEvaluated}
    \detTrafoEvaluated
    }
    + \reciprocal{\pecletNumber} \IntegralDomain{\transpose{(\invTranspJacTrafoEvaluated  \gradient{\testSpaceSturmApproachVolumePercentageSolid})} (\invTranspJacTrafoEvaluated  \gradient{\ansatzSpaceSturmApproachVolumePercentageSolid}) \detTrafoEvaluated} \nonumber \\ 
      &+\IntegralSurfacePartDomainPlain{\transpose{(\ansatzSpaceSturmApproachVelocity- \velocityFluidInflow)} \, \testSpaceSturmApproachVelocityFluidInflowCondition \detTrafoBdryEvaluated}{\surfaceCompDomainInflow}     
    +\IntegralSurfacePartDomainPlain{\transpose{(\ansatzSpaceSturmApproachVelocitySolid - \velocitySolidInflow)} \testSpaceSturmApproachVelocitySolidInflowCondition \detTrafoBdryEvaluated}{\surfaceCompDomainInflow} 
    +\IntegralSurfacePartDomainPlain{(\ansatzSpaceSturmApproachVolumePercentageSolid - \volumePercentageSolidInflow)\testSpaceSturmApproachVolumePercentageSolidInflowCondition \detTrafoBdryEvaluated}{\surfaceCompDomainInflow}. \nonumber
\end{align}
\endgroup

\Cref{eq:ShapeDerivativeEulerEulerModel}
is then obtained from \cref{eq:ComputationOfEulerianSemiderivative}
using the identities \cref{eq:DifferentiationOfTransformationQuantities,eq:TotalDerivativeOuterNormalOnTransformedDomainWithPushForward},
where the last three integrals in \cref{eq:TransformedLagrangianEulerEulerModel} vanish
due to
$\detTrafoBdryPrimeAtZero = \tangentialDivergence{\boundaryVelocity}=0$
on $\surfaceCompDomainInflow$.
\end{pf}

%% file: GradientDescendMethod_and_Discretization.tex
This section is devoted to the description of the numerical 
framework for the gradient-based treatment of the shape optimization problem \cref{eq:OptimizationProblem}.
In \cref{sec:GradientProjection} we describe the projection of the shape derivative \cref{eq:ShapeDerivativeEulerEulerModel} 
and the mesh deformation procedure,
and in \cref{sec:DiscretizationAndSolvers} we comment on the discretization of 
the partial differential equations
with finite elements, as well as the non-linear and linear solvers.

\subsection{Gradient projection}
\label{sec:GradientProjection}

In order to obtain smooth mesh deformations from the shape derivative \cref{eq:ShapeDerivativeEulerEulerModel}, 
we use the approach of solving the linear elasticity equations with the volume and surface parts of the shape gradient acting as body forces and surface tractions, see e.g. \cite{dwight2009robust,schulz2016computational},
in conjunction with a recently proposed correction procedure for the volume mesh deformation \cite{etling2018first}.
Given the space of deformations
\begin{equation*}
 \hilbertSpaceDeformation
 := \left\{ W \in \HOneDomainVect \,\middle|\, \evaluateAtTwo{W}{\surfaceCompDomain \setminus \surfaceCompDomainDeformable} = 0 \right\},
\end{equation*}
we identify the shape derivative 
\cref{eq:ShapeDerivativeEulerEulerModel} 
with an element $\gradientInHilbertSpace \in \hilbertSpaceDeformation$ 
by solving 
\begin{equation}
\label{eq:ShapeGradientProjection}
 \dependence{\bilinearFormProjection}{\gradientInHilbertSpace,\boundaryVelocity}
 = \denoteEulerianDerivativeFunctionalDomain{\costFunctionalLetter}{\boundaryVelocity}
 \quad \forall \, \boundaryVelocity \in \hilbertSpaceDeformation 
\end{equation}
with an elliptic bilinear form $\bilinearFormProjection$.
For this purpose we introduce the strain and stress tensors
\begin{equation*}
 \dependence{\strainTensor}{\gradientInHilbertSpace} 
 := \lameParameterOne \trace{\dependence{\stressTensor}{\gradientInHilbertSpace}}\identityMatrix
   + 2 \lameParameterTwo \dependence{\stressTensor}{\gradientInHilbertSpace},
   \quad
\dependence{\stressTensor}{\gradientInHilbertSpace}
:= \frac{1}{2} (\gradient{\gradientInHilbertSpace} + \transpose{\gradient{\gradientInHilbertSpace}}), 
\end{equation*}
with the first and second Lam{\'e} parameters $\lameParameterOne$ and $\lameParameterTwo$.
Proceeding in a similar manner as \cite{schulz2016computational},
we set $\PoissonRation= \lameParameterOne = 0$
and 
$\lameParameterTwo:=\sqrt{\lameParameterTwoPDE}$ 
with the solution $\lameParameterTwoPDE$ of the weak form of 
\begin{align}
\begin{split}
  \label{eq:ExtensionLameParameterTwo}
 \laplace{\lameParameterTwoPDE} &= 0 \textup{ in } \domain, \\
 \lameParameterTwoPDE &= \lameParameterTwoPDEMax \textup{ on } \surfaceCompDomainDeformable, \\
 \lameParameterTwoPDE &= \lameParameterTwoPDEMin \textup{ on } \surfaceCompDomain \setminus \surfaceCompDomainDeformable,
\end{split}
\end{align}
for $\lameParameterTwoPDEMin=1$ and $\lameParameterTwoPDEMax=100$,
and set 
\begin{equation}
\label{eq:DefinitionLinearElasticityInnerProduct}
\dependence{\bilinearFormProjection}{\gradientInHilbertSpace,\boundaryVelocity}
:=
   \IntegralDomain{\dependence{\strainTensor}{\gradientInHilbertSpace} : \dependence{\stressTensor}{\boundaryVelocity}}.
\end{equation}

In addition to the projection step \cref{eq:ShapeGradientProjection},
we use a recently proposed gradient correction procedure \cite{etling2018first},
which has been shown to successfully remove spurious components in the discretized gradient.
Upon introducing 
\begin{align*}
\defineFunctionLong{\normalForceOperator}
		   {\LTwoBoundary \times \HOneDomainVect}
		   {\R}
		   {(\normalForceLagrangeParameter,\boundaryVelocity)}
		   {
		      \IntegralSurfaceDomainPlain{\normalForceLagrangeParameter (\transpose{\boundaryVelocity }\outerNormal)},
		    }
\end{align*}
this method introduces the additional saddle-point problem
\begin{align}
 \begin{split}
  \label{eq:SaddlePointProblemShapeGradientCorrection}
  \textup{Find } (\normalForceLagrangeParameter,\boundaryVelocityCorrection) \in \LTwoBoundary &\times \HOneDomainVect \textup{, s.t. } \\
  \dependence{\normalForceOperator}{\normalForceLagrangeParameterTestFunction,\boundaryVelocityCorrection}
    &= 0 \qquad\qquad\,\,  \forall \normalForceLagrangeParameterTestFunction \in \LTwoBoundary, \\
   \dependence{\normalForceOperator}{\normalForceLagrangeParameter,\boundaryVelocity}
   + \dependence{\bilinearFormProjection}{\boundaryVelocityCorrection,\boundaryVelocity}
    &= \denoteEulerianDerivativeFunctionalDomain{\costFunctionalLetter}{\boundaryVelocity}
 \quad \forall \, \boundaryVelocity \in \hilbertSpaceDeformation.
  \end{split}
\end{align}

Given the solutions 
$\gradientInHilbertSpaceDiscretized$ and 
$(\normalForceLagrangeParameterDiscretized,\boundaryVelocityCorrectionDiscretized)$
of the 
Galerkin finite element formulations of  \cref{eq:ShapeGradientProjection,eq:SaddlePointProblemShapeGradientCorrection}
with linear Lagrange elements, 
the components of the discrete restricted shape gradient 
\begin{equation}
\label{eq:DefinitionGradientRestrictedMeshDeformation}
\gradientInHilbertSpaceCorrectedDiscretized:=\gradientInHilbertSpaceDiscretized-\boundaryVelocityCorrectionDiscretized \in \hilbertSpaceDeformationFiniteDimensional
\end{equation}
contain the deformation directions at the $\expDiscrHilbertSpaceDef\in\N$ mesh vertices 
$\meshVertices \in \hilbertSpaceDeformationFiniteDimensional$.
Note, that due to the first equation in \cref{eq:SaddlePointProblemShapeGradientCorrection}
the correction variable $\boundaryVelocityCorrection$ contains only tangential components at $\surfaceCompDomain$,
and that $\gradientInHilbertSpaceCorrectedDiscretized$ and $\gradientInHilbertSpaceDiscretized$ therefore induce the same shape changes up to the discretization error introduced by the approximation of this saddle point problem.

With \cref{eq:DefinitionGradientRestrictedMeshDeformation} and the discretized formulation of 
the forward and adjoint PDEs as well as  \cref{eq:ExtensionLameParameterTwo,eq:ShapeGradientProjection,eq:SaddlePointProblemShapeGradientCorrection},
we iteratively deform the current mesh according to \Cref{alg:GradientDescentMethod}.
Here the step size $\denoteIteration{\stepsizeBacktracking}>0$
is required to fulfill the Armijo condition \cite{wright1999numerical} 
in order to obtain a sufficient decrease of $\costFunctionalLetter$, 
and 
we additionally impose \cite{etling2018first} 
\begin{equation}
  \label{eq:StepsizeConditionEtling}
 \frac{1}{2} \le \det{(\identityMatrix + \denoteIteration{\stepsizeBacktracking} \jacobian{\denoteIteration{\gradientInHilbertSpaceCorrectedDiscretized}})} \le 2,
 \qquad \norm{\denoteIteration{\stepsizeBacktracking} \jacobian{\denoteIteration{\gradientInHilbertSpaceCorrectedDiscretized}}}_{F} \le 0.3,
\end{equation}
with the Frobenius norm $\norm{M}_{F}$ of the matrix $M$.
These restrictions intend to restrict the maximal mesh volume and angle changes
in the deformation step and therefore to avoid inverted mesh elements.

\begin{alg}[Gradient descent method]
 \label{alg:GradientDescentMethod}
 \begin{algorithmic}
   \STATE{$ $}\\ 
   \STATE{Create mesh $\meshVertices_{1}$ on the initial domain $\domain_1$}\\
   \STATE{Set $converged \leftarrow false$}\\
  \REPEAT
   \STATE{Compute state and adjoint variables from \eqref{eq:WeakForwardPDESystem} and \eqref{eq:AdjointEquationFromShapeLagrangian}}\\
   \STATE{Compute restricted shape gradient $\denoteIteration{\gradientInHilbertSpaceCorrectedDiscretized}$ from \cref{eq:ExtensionLameParameterTwo,eq:ShapeGradientProjection,eq:SaddlePointProblemShapeGradientCorrection}}\\ 
    \IF{Step size $\denoteIteration{\stepsizeBacktracking}$ yields decrease in the cost function}
    \STATE{$\denoteNextIteration{\meshVertices} \leftarrow
    \denoteIteration{\meshVertices} - \denoteIteration{\stepsizeBacktracking} \denoteIteration{\gradientInHilbertSpaceCorrectedDiscretized}    
    $}\\
    \ELSE
      \STATE{$converged \leftarrow true$}
    \ENDIF
    \UNTIL{$converged$}
\end{algorithmic}
\end{alg}

\begin{rmk}
 We note, that there exist more sophisticated higher-order methods for shape optimization problems based on the analytical shape Hessian or approximations thereof; 
 see  \cite{eppler2005regularized,etling2018first,schmidt2009impulse,schulz2016computational,schulz2015structured} for some examples.
 However, since \Cref{alg:GradientDescentMethod} 
 yielded satisfactory results for the test case considered in \cref{sec:NumericalResults}, 
 we did not further pursue one of these approaches.
\end{rmk}

\subsection{Discretization and iterative solvers}
\label{sec:DiscretizationAndSolvers}

For the implementation of \Cref{alg:GradientDescentMethod}
we utilize the assembling and solution capabilities of COMSOL Multiphysics 5.3a \cite{comsol}.
In order to keep the computational effort for one iteration of \Cref{alg:GradientDescentMethod} on a given mesh as small as possible,
we choose linear Lagrange elements for the discretization of the state and adjoint variables as well as 
$\gradientInHilbertSpaceDiscretized$, 
$\normalForceLagrangeParameterDiscretized$ 
and 
$\boundaryVelocityCorrectionDiscretized$.

It is well known, that equal-order Lagrange elements for fluid velocity and pressure do not fulfill the LBB-condition \cite{ErnGuermond2004} and stabilization techniques have to be considered. 
For the forward and adjoint fluid equations we apply the 
Streamline-Upwind-Petrov-Galerkin (SUPG) 
and  
Pressure-Stabilized-Petrov-Galerkin (PSPG) methods,
as well as a 
least-squares penalization of the incompressibility condition (LSIC)
 with the parameters proposed in  \cite{peterson2018overview}.
The forward and adjoint particle velocity equations are stabilized through
a combination of the SUPG-method \cite{bourgault1999finite} 
and a viscosity ramping strategy for the diffusion constant up to $\reciprocal{\smoothingVelocitySolidNumber}=10^{-4}$.
In contrast to the Navier-Stokes and particle velocity equations, 
we use the predefined COMSOL module \textit{Chemical Species Transport Interface}
for the forward and adjoint transport equations.
This module includes SUPG and crosswind-diffusion stabilization \cite{brooks1982streamline,do1991feedback},
the later of which introduces non-linear terms to these linear equations.

The non-linear forward equations as well as the adjoint transport equation
are solved with a damped Newton method,
where a restarted GMRES solver with Krylov-dimension 50 and an algebraic multigrid preconditioner is used for the linear sub-problems,
and the same linear solver is also used for the remaining adjoint equations.
Due to the ellipticity of the differential operators in \cref{eq:ShapeGradientProjection} and \cref{eq:ExtensionLameParameterTwo},
we use the CG method in conjunction with a symmetric overrelaxed Gauss-Seidel preconditioner for the corresponding discretized problems. 
Owing to its saddle-point structure, 
this preconditioning technique can not be applied 
for the volume mesh correction step \cref{eq:SaddlePointProblemShapeGradientCorrection}, 
where the symmetric system matrix of the discretized problem is not positive definite.
Despite this restriction, 
the CG method without preconditioner showed satisfactory convergence behavior for the test cases considered in \cref{sec:NumericalResults}, 
which is why we preferred it over, e.g., an approach based on the Schur complement \cite{zhang2006schur}.

%% file: NumericalResults.tex
This section is devoted to the numerical validation 
of the Eulerian particle model 
and the shape derivative \cref{eq:ShapeDerivativeEulerEulerModel} 
of \cref{eq:CostFunctionReduced}.
To this end, we compare the predicted particle impact rates with respect to the Stokes number for a $90^{\circ}$ bend to experimental data and reference simulations in \cref{sec:ValidationImpactRates}.
In \cref{sec:ShapeOptimizationOfPipeBend} we then 
apply \Cref{alg:GradientDescentMethod} 
to minimize the erosion rate caused by a monodisperse particle jet 
for this test case and compare the initial and optimized geometry for a wider range of Stokes numbers.

%
%
\subsection{Validation of computed impact rates}
\label{sec:ValidationImpactRates}
\begin{figure}[htbp!]
  \begin{minipage}[t]{0.32\textwidth}
    \centering
     \includegraphics[height=7cm]{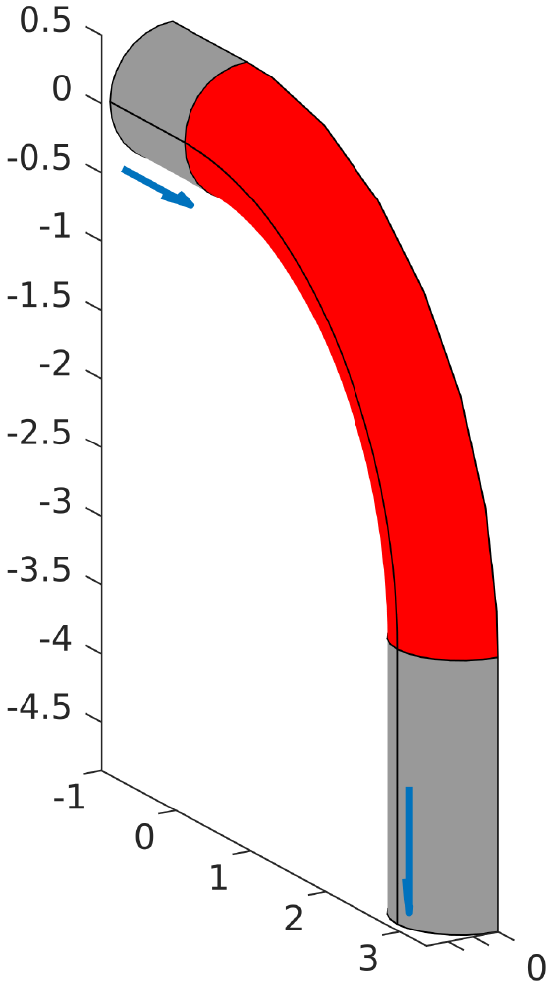}
     \caption{Initial geometry and deformable boundary part}
     \label{fig:InitialGeometry}
  \end{minipage}
  \hfill
  \begin{minipage}[t]{0.64\textwidth}
      \centering
	\begin{tikzpicture}[scale=1.0]
	\begin{axis}[
	    ymin=0,
	    ymax=1,
	    legend pos=south east,
	    legend cell align={left},
	    xlabel=$\stokesNumber$,
	    ylabel=$\depositionFraction$
	  ]
	\addplot [mark=*, mark options={solid}]
	    table[]{Graphics/tikz/data/DepositionRateOutflowToInflow_Initial.tsv}; 
	    \addlegendentry{Eulerian model \cref{eq:DepositionFraction}}    
	\addplot [mark=square*]
	  table [x=XPui1987, y=YPui1987,           col sep=comma] {Data_Deposition_literature/wpd_datasets.csv};
	  \addlegendentry{Pui et al. \cite{pui1987experimental}}
	\addplot [mark=triangle*]
	  table [x=XVasquez, y=YVasquez,           col sep=comma] {Data_Deposition_literature/wpd_datasets.csv};
	  \addlegendentry{Vasquez et al. \cite{vasquez2015analysis}}
	\addplot [mark=none,densely dashdotted]
	  table [x=XPilou2011, y=YPilou2011,       col sep=comma] {Data_Deposition_literature/wpd_datasets.csv};
	  \addlegendentry{Pilou et al. \cite{pilou2011inertial}}
	\addplot [mark=none,loosely dashed]
	  table [x=XTsai_Pui1990, y=YTsai_Pui1990, col sep=comma] {Data_Deposition_literature/wpd_datasets.csv};
	  \addlegendentry{Tsai \& Pui   \cite{tsai1990numerical}}
	\addplot [mark=none,dotted]
	  table [x=XBreuer2006, y=YBreuer2006,     col sep=comma] {Data_Deposition_literature/wpd_datasets.csv};
	  \addlegendentry{Breuer et al. \cite{breuer2006prediction}}
	\end{axis}
      \end{tikzpicture}
    \caption{Comparison of calculated impact rates with numerical and experimental studies}
    \label{fig:DepositionFractionsComparisonToLiterature}
  \end{minipage}
\end{figure}

In order to verify the accuracy of the non-conservative formulation of the Eulerian particle model,
we compute the impact rates 
for various Stokes numbers by changing the particle diameters
for the bend depicted in \Cref{fig:InitialGeometry},
where due to the symmetry of the problem only one half of the geometry has to be considered.
The complete set of parameters describing the geometry 
based on \Cref{fig:PlotBendedPipe}
as well as the flow and particle parameters is given in \Cref{tab:ParametersCollectionEfficiency}. The parameters are chosen in accordance with the experimental test case considered 
\cite{pui1987experimental} 
and the numerical studies
\cite{breuer2006prediction,pilou2011inertial,tsai1990numerical,vasquez2015analysis} thereof.

\begin{table}[htbp!]
\centering
\caption{Fluid, particle and stabilization parameters (in SI-units)}
\label{tab:ParametersCollectionEfficiency}
\begin{tabular}{l|l|l}
 \toprule
 Variable & Value & Description \\
 \midrule
  $\dTube$ 
    & $3.95 \cdot 10^{-3}$ & inlet and outlet diameter\\
 $\rCurvatureBend$ 
    & $1.13 \cdot 10^{-2}$ & radius of curvature\\
 $\rCurvatureRatio$ 
    & $5.7$ & curvature ratio\\
 $\densityFluid$ 
    & $1.18$ & fluid density\\
 $\viscosityFluid$ 
    & $1.85 \cdot 10^{-5}$ & dynamic fluid viscosity\\
 $\velocityFluidMeanInflow$ 
    & $3.86$ & mean fluid velocity on $\surfaceCompDomainInflow$\\    
 $\densitySolid$ 
    & $895$ & density of particles\\
 $\diamEq$ 
    & $5-16 \cdot 10^{-6}$ & particle diameter\\
 $\stokesNumber$ 
    & $0.13-1.34$ & Stokes number\\    
  $\froudeNumber$ 
    & $27.7$ & Froude number\\   
 $\reynoldsNumber$ 
    & $1000$ & Reynolds number\\
 $\deanNumber$
    & 419 & Dean number \\
 $\smoothingVelocitySolidNumber$ 
    & $10^{4}$ & artificial viscosity constant\\
 $\pecletNumber$ 
    & $10^{8}$ & Peclet number\\    
    \bottomrule
\end{tabular}
\end{table}

The question of an adequate mesh size for this test case has been investigated in \cite{vasquez2015analysis}, 
where the authors report, 
that a resolution beyond approximately $5 \times 10^5$ tetrahedral elements for the full bend 
including a boundary layer refinement at $\surfaceCompDomainWall$ 
did not change the predicted impact rates. 
In order to obtain a satisfying spatial resolution of the derivatives of the erosion rate $\functionalErosionLetter$ on $\surfaceCompDomainWall$,
we use a slightly finer mesh, which consists of $1.28\times10^6$ tetrahedral elements and $\expDiscrHilbertSpaceDef=3.2\times 10^5$ vertices with a prismatic boundary layer refinement at $\surfaceCompDomainWall$ for the half bend from \Cref{fig:InitialGeometry}.

Given the solution $\solutionSturmApproachLetter \in \spaceE$ of \cref{eq:WeakForwardPDESystem}, 
the predicted particle impact rates 
\begin{align}
  \label{eq:DepositionFraction}
 \depositionFraction 
   :=  1- 
   \IntegralSurfacePartDomainPlain{\solutionSturmApproachVolumePercentageSolid (\transpose{\solutionSturmApproachVelocitySolid} \outerNormal)}{\surfaceCompDomainOutflow}
    \bigg/ 
 \IntegralSurfacePartDomainPlain{\solutionSturmApproachVolumePercentageSolid (\transpose{\solutionSturmApproachVelocitySolid} (-\outerNormal))}{\surfaceCompDomainInflow} 
\end{align}
are depicted in \Cref{fig:DepositionFractionsComparisonToLiterature}.
While for $\stokesNumber>0.43$ very good agreement with the reference studies can be observed,
the Eulerian model seems to slightly over-predict $\depositionFraction$ for $\stokesNumber \le 0.43$.
However, the deviation of the values corresponding to these small diameters 
from the highest impact rates reported in \cite{breuer2006prediction} 
is within the range of deviations among the reference studies,
so that we proceed with the treatment of the optimization problem introduced in \cref{sec:OptimizationProblem}.

%
%
\subsection{Shape optimization of the bend}
\label{sec:ShapeOptimizationOfPipeBend}

In this section we apply \Cref{alg:GradientDescentMethod}
to the initial geometry from \Cref{fig:InitialGeometry}
for the numerical solution of the optimization problem \cref{eq:OptimizationProblem}.
In order to demonstrate our approach, we
use the model parameters given in \Cref{tab:ParametersCollectionEfficiency}
for the particle species with $\stokesNumber=0.33$.
Since we want to keep the influence of the curvature regularization on the cost functional \cref{eq:CostFunctionalEulerEulerModel} and thus on the optimized shape as small as possible, 
we choose
$\constantRegularizationOne = 0.01 \IntegralSurfacePartDomainPlain{\dependence{\costFunctionalBoundaryFunction}{\volumePercentageSolid,\velocitySolid,\outerNormal}}{\surfaceCompDomainWall} /  \IntegralSurfacePartDomainPlain{\frac{1}{2} \, \meanCurvature^2}{\surfaceCompDomainDeformable}$,
where the integrals are evaluated on the initial geometry.
Additionally we use the parameters  
$\exponentVelocityDependencyOka=2.36$, $\hardnessOkaModel=2$,
$\exponentBrittleOkaModel=\exponentBrittleOkaModelDuctileCaseVal$
and
$\exponentDuctileOkaModel=\exponentDuctileOkaModelDuctileCaseVal$
for the erosion rate \cref{eq:ErosionRateOKA}, 
which are taken from \cite{oka1997impact} to model the material characteristics of stainless steel.
With these parameters, the local erosion rate is smallest for almost tangential impacts.

Inspired by the Hadamard structure theorem 
\cite{Sokolowski1992introduction},
we consider the norm
\begin{equation}
\label{eq:GradientNorm}
 \norm{\gradientInHilbertSpaceCorrected}_{\surfaceCompDomainDeformable}
 := \sqrt{
      \IntegralSurfacePartDomainPlain{(\gradientInHilbertSpaceCorrected \cdot \outerNormal)^2}{\surfaceCompDomainDeformable}
    }
\end{equation}
as a measure for shape changes induced by 
$\gradientInHilbertSpaceCorrected \in \hilbertSpaceDeformation$
and iteratively deform the current geometry according to \Cref{alg:GradientDescentMethod}
until no more decrease in the objective can be obtained. 
To ensure, that this only happens close to a stationary point of the optimization problem \cref{eq:OptimizationProblem}, we track the relative decrease in the objective as well as the gradient norms 
$\denoteIteration{\gradientInHilbertSpaceCorrectedDiscretized}  
\in \hilbertSpaceDeformationFiniteDimensional$ 
with the discrete counterpart of \cref{eq:GradientNorm} on the corresponding geometry. 
Since we only consider piece-wise linear elements for the representation of the discrete shape derivative, 
an approximation of the integrals in \cref{eq:GradientNorm} with the trapezoidal rule is sufficiently accurate.

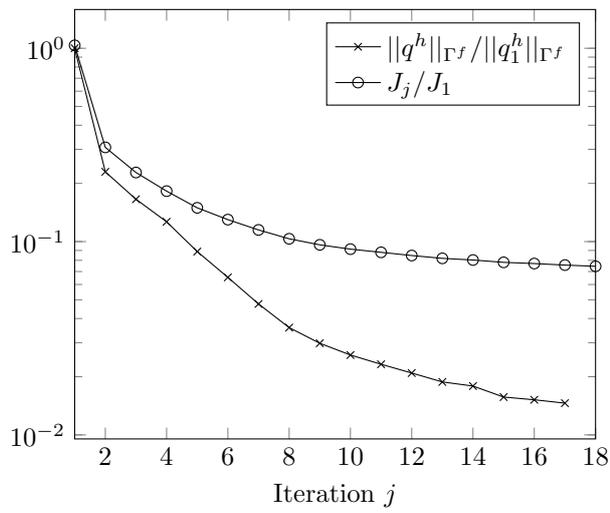
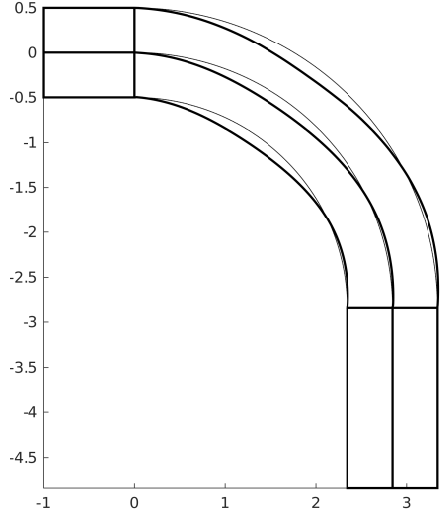
\begin{figure}[]
  \begin{subfigure}[t]{0.48\textwidth}
   \centering
	  \begin{tikzpicture}[scale=1.0]
	  \begin{axis}[
	      xmin=1,
	      xmax=18,
	      ymode=log,
	      legend pos=north east,
	      legend cell align={left},
	      xlabel=$\textup{Iteration } \iterBFGS $
	    ]
	
	  \addplot [mark=x, mark options={solid}]
	      table[]{Graphics/tikz/data/Relative_d_dir_norms_3BFGS.tsv}; 
	      \addlegendentry{$\norm{\lbfgsDescentDirectionDiscrete}_{\surfaceCompDomainDeformable} / \norm{\lbfgsDescentDirectionDiscrete_1}_{\surfaceCompDomainDeformable} $}
	      
	  \addplot [mark=o, mark options={solid}]
	      table[]{Graphics/tikz/data/Relative_cost_functional_3BFGS.tsv}; 
	      \addlegendentry{$\denoteIteration{\costFunctionalLetter} / \costFunctionalLetter_1$}
	      \end{axis}
	  \end{tikzpicture}
    \caption{Convergence history of the optimization algorithm for $\stokesNumber=0.33$}
    \label{fig:ConvergenceHistory}
  \end{subfigure}
  ~
  \begin{subfigure}[t]{0.48\textwidth}
    \centering
    \includegraphics[trim={9.5em 9.5em 9.0em 9.0em}, clip, scale=0.56]{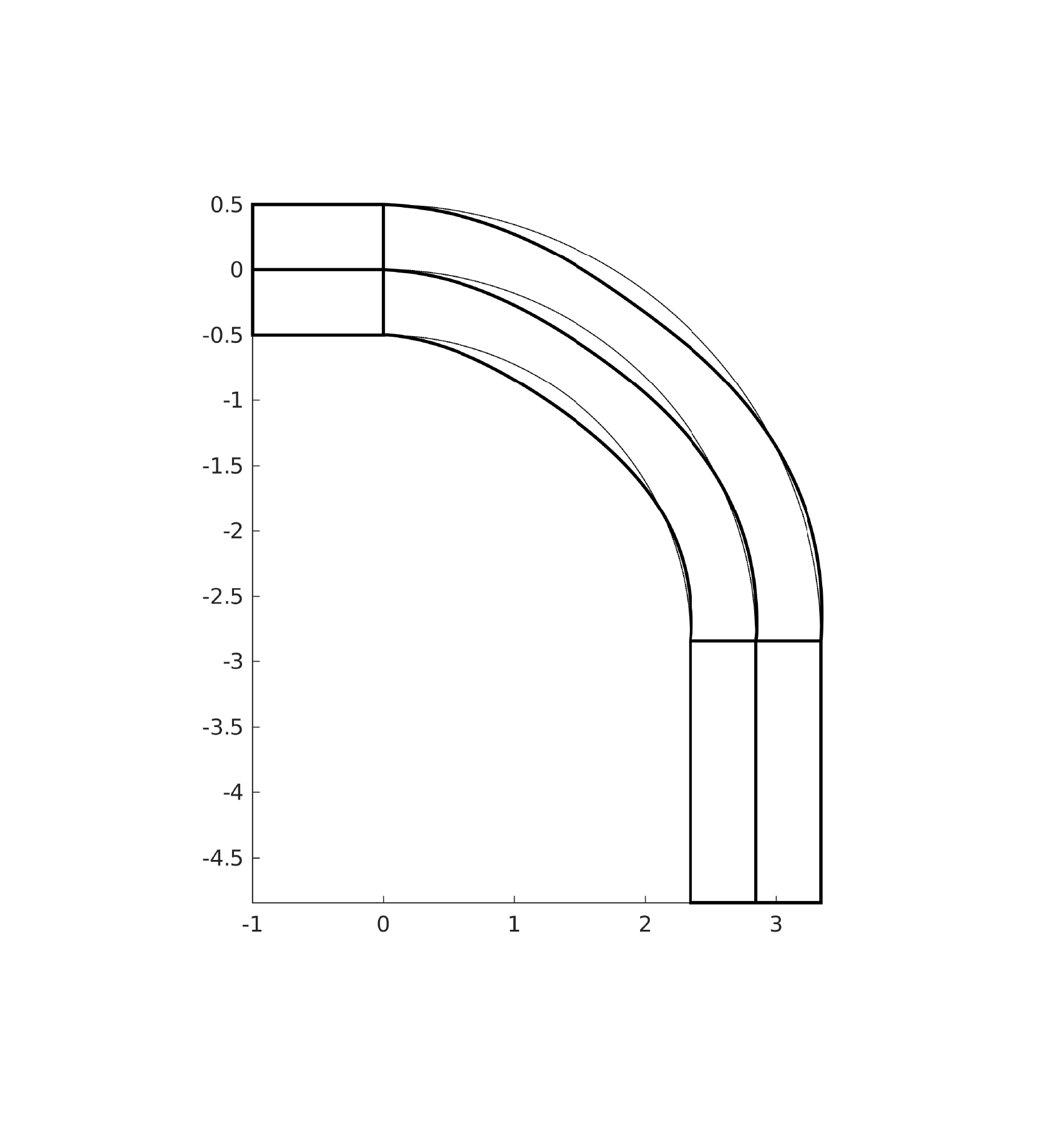}
    \subcaption{Optimized (thick lines) and initial geometry (thin lines)}
    \label{fig:InitialAndOptimizedGeometry}
  \end{subfigure}
 \caption{Application of \Cref{alg:GradientDescentMethod} to the initial geometry from \Cref{fig:InitialGeometry} }
 \label{fig:ConvergenceHistoryAndInitialOptimizedGeometry}
\end{figure}

\begin{figure}[htbp!]
  \begin{subfigure}[t]{0.48\textwidth}
	\centering
	\includegraphics[height=8cm]{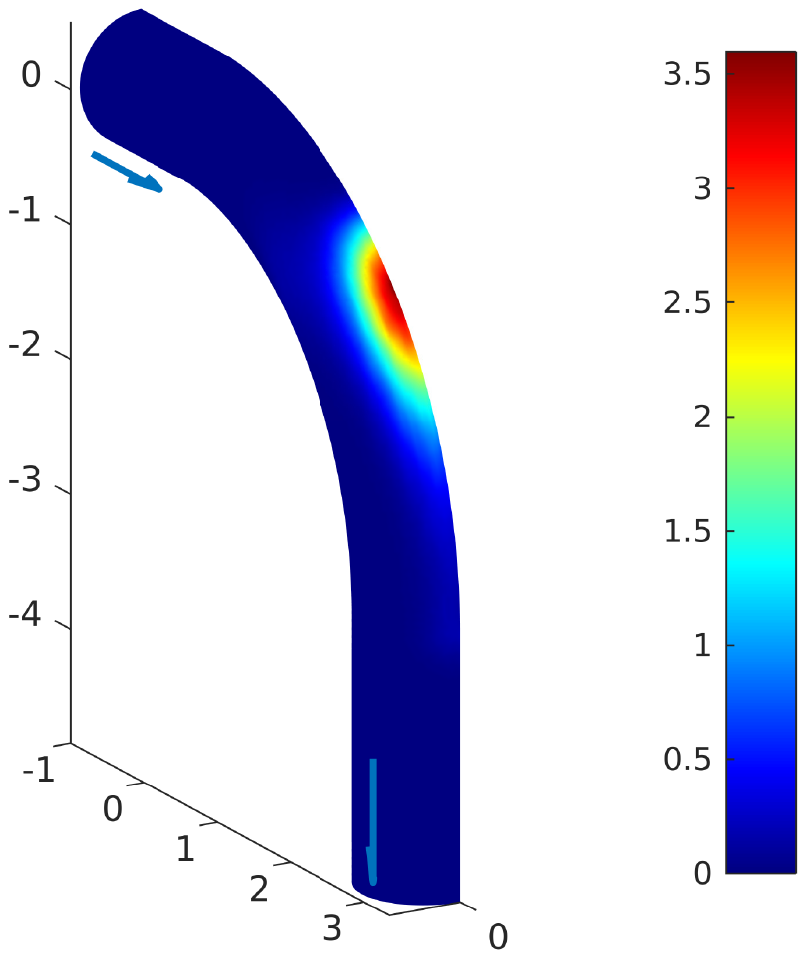}
	\subcaption{Erosion rate $\functionalErosionLetter$ on the initial bend}
	\label{fig:ErosionRateInitial}
    \end{subfigure}
    ~
     \begin{subfigure}[t]{0.48\textwidth}
	\centering
	\includegraphics[height=8cm]{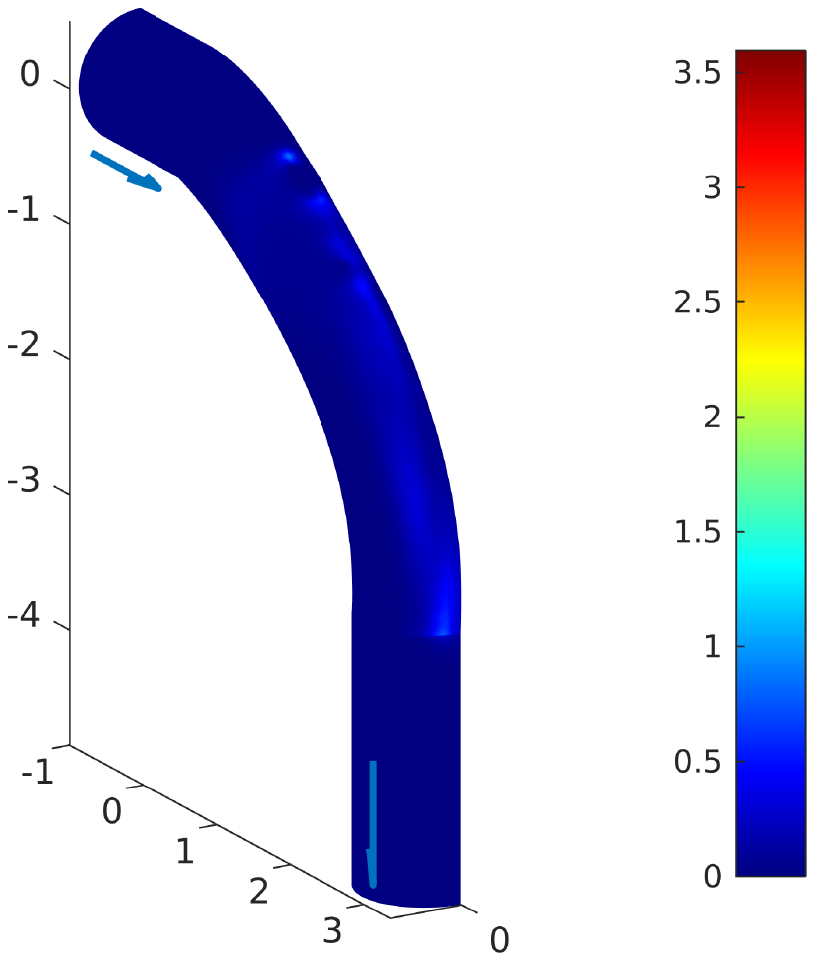}
	\subcaption{Erosion rate $\functionalErosionLetter$ on the optimized bend}
	\label{fig:ErosionRateOptimized}
    \end{subfigure}
    \caption{Erosion rates on $\domain_1$ and $\domain_{18}$}
    \label{fig:ErosionRatesOnInitialAndOptimizedGeometries}
\end{figure}

\Cref{fig:ConvergenceHistory} shows the decrease in the objective and the gradient norms  \cref{eq:GradientNorm} during the optimization. 
Since most of the decrease in the objective is achieved in the first half of the iterations,
and since the relative gradient norms fall even beyond that point, 
we deduce that a locally optimal shape has been obtained.
This geometry is depicted in \Cref{fig:InitialAndOptimizedGeometry} together with the initial bend,
and in \Cref{fig:ErosionRatesOnInitialAndOptimizedGeometries} the erosion rates on both geometries are compared. 
Since we incorporated the squared erosion rate in the definition of the cost functional \cref{eq:CostFunctionalEulerEulerModel},
the maximal erosion rate due to the impact of particles with $\stokesNumber=0.33$
can be decreased 
by $76\%$ with respect to the initial geometry.

\begin{figure}[htbp!]
  \begin{subfigure}[t]{0.32\textwidth}
	\centering
	\includegraphics[height=5.6cm]{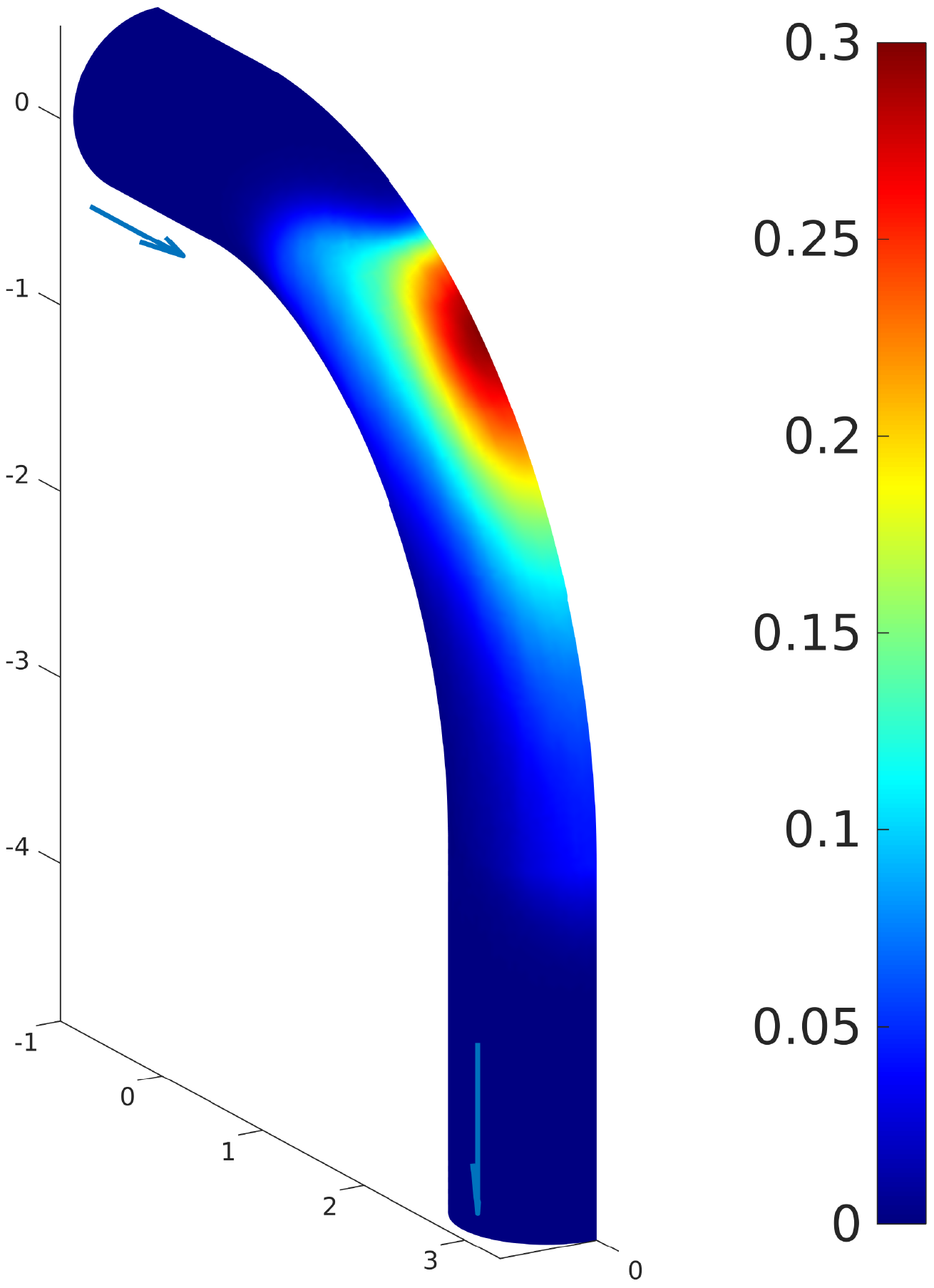}
	\subcaption{Impact rate $\volumePercentageSolid(\transpose{\velocitySolid}\outerNormal)$}
	\label{fig:ImpactRateInitial}
    \end{subfigure}
    ~
     \begin{subfigure}[t]{0.32\textwidth}
	\centering
	\includegraphics[height=5.6cm]{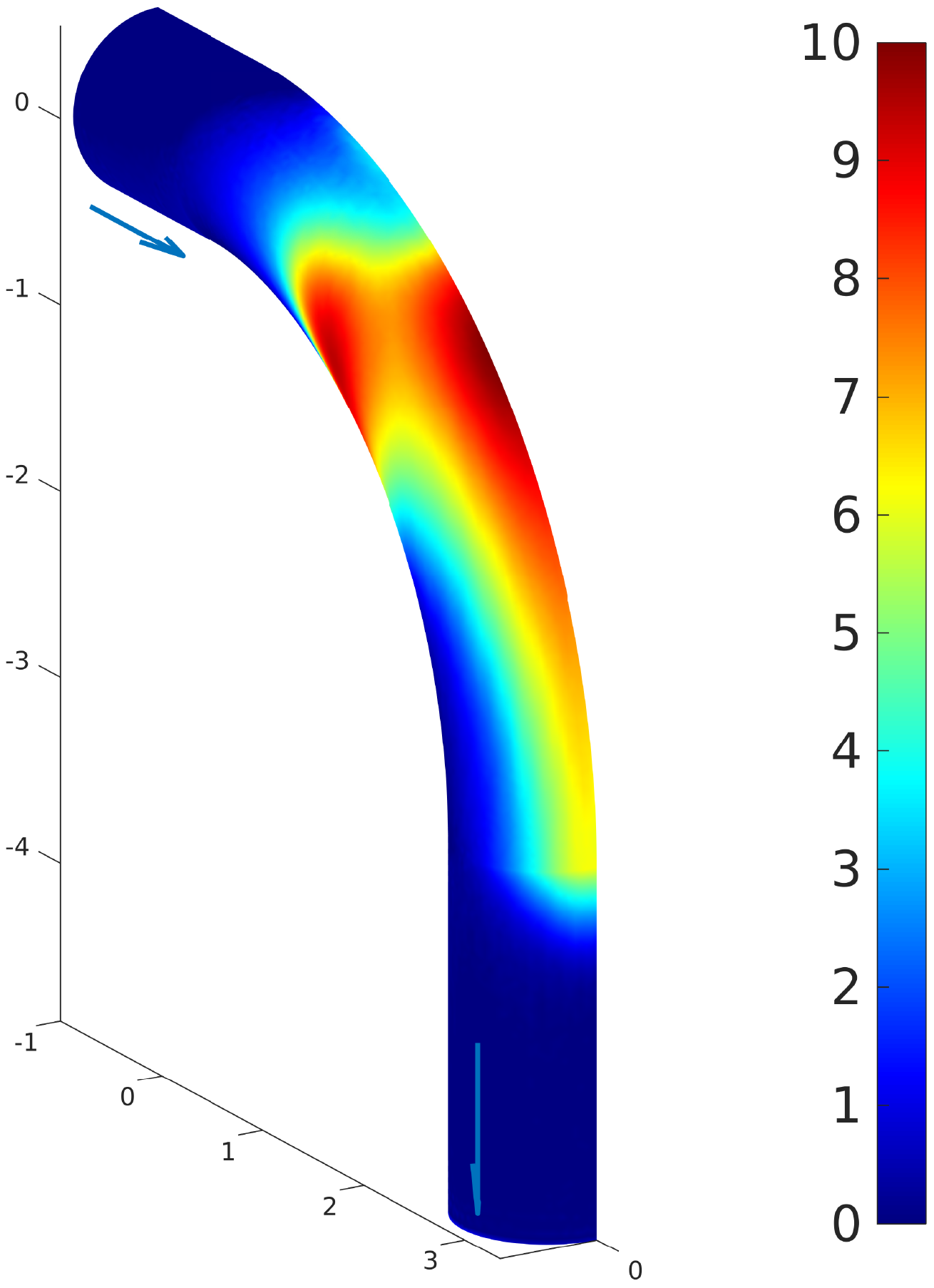}
	\subcaption{Impact angle $\dependence{\angleErosion}{\velocitySolid,\outerNormal}$ (in $^{\circ}$)}
      \label{fig:ImpactAnglesInitial}
    \end{subfigure}
    ~
     \begin{subfigure}[t]{0.32\textwidth}
	\centering
	\includegraphics[height=5.6cm]{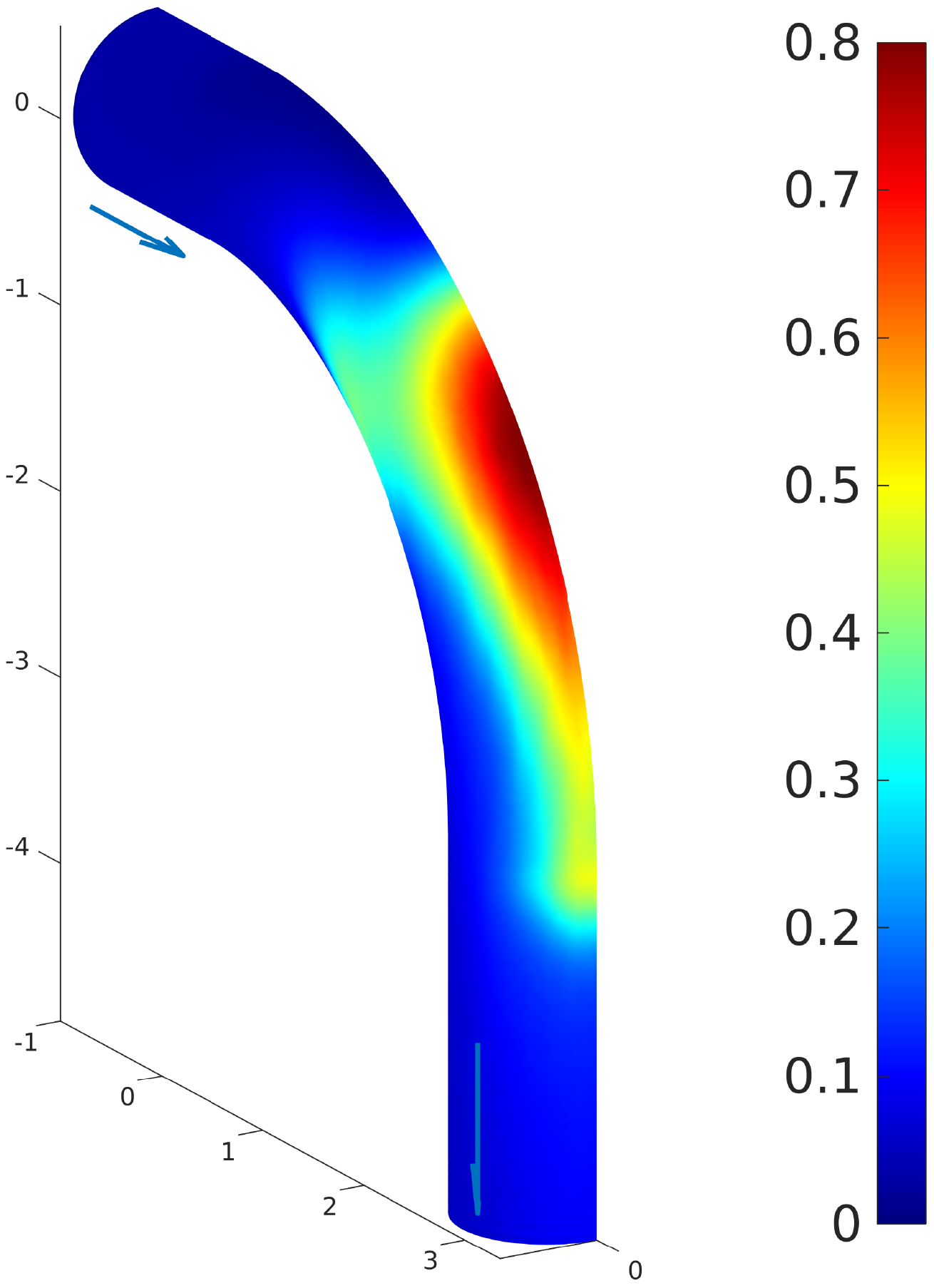}
	\subcaption{Impact velocity $\normEuclidean{\velocitySolid}$}
	\label{fig:ImpactVelocitiesInitial}
    \end{subfigure}
    \caption{Impact rates, velocities and angles for the initial geometry}
    \label{fig:ComponentsOfErosionRateInitialGeometry}
\end{figure}

\begin{figure}[htbp!]
  \begin{subfigure}[t]{0.32\textwidth}
	\centering
	\includegraphics[height=5.6cm]{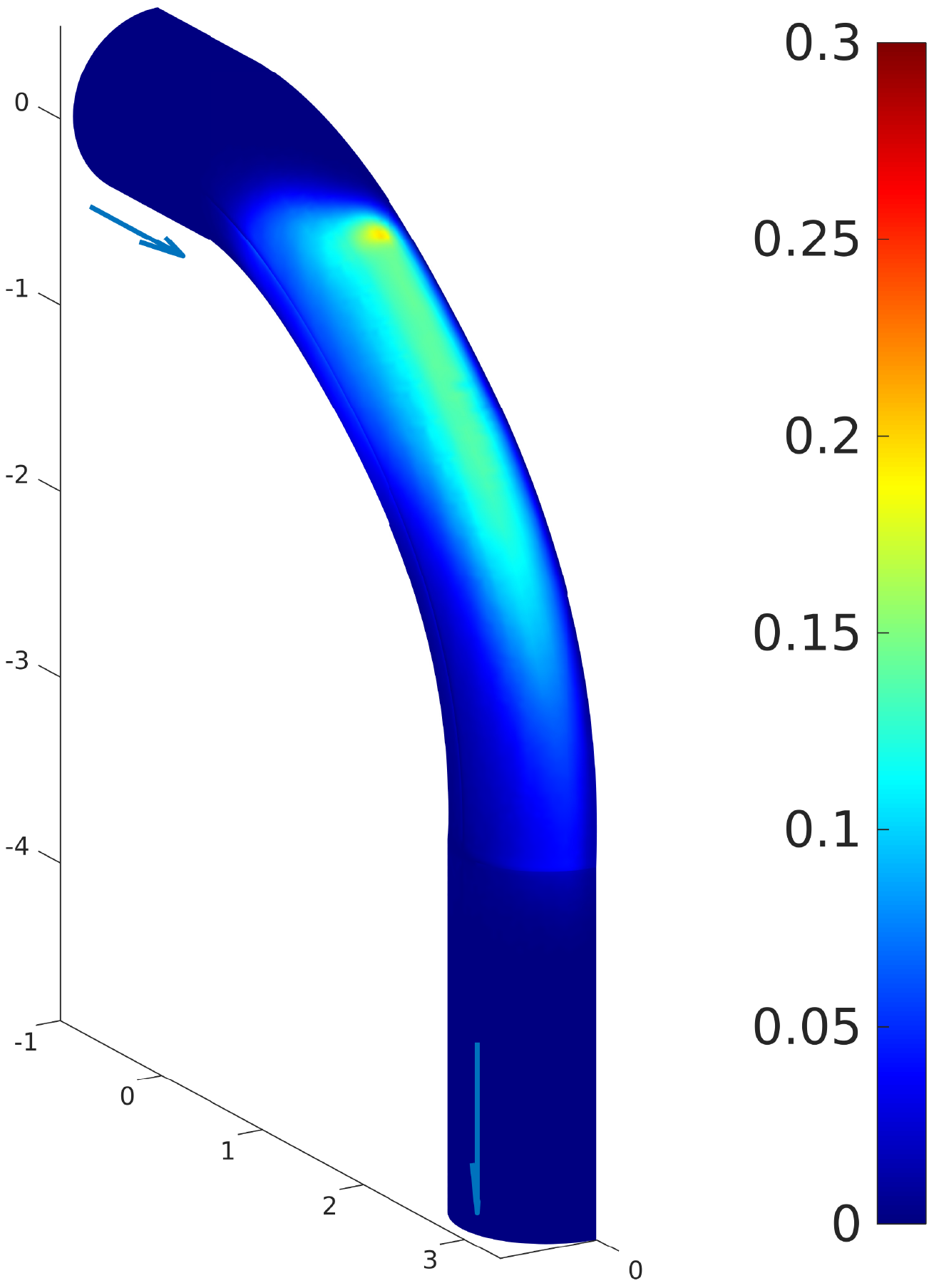}
	\subcaption{Impact rate $\volumePercentageSolid(\transpose{\velocitySolid}\outerNormal)$}
	\label{fig:ImpactRateOptimized}
    \end{subfigure}
    ~
     \begin{subfigure}[t]{0.32\textwidth}
	\centering
	\includegraphics[height=5.6cm]{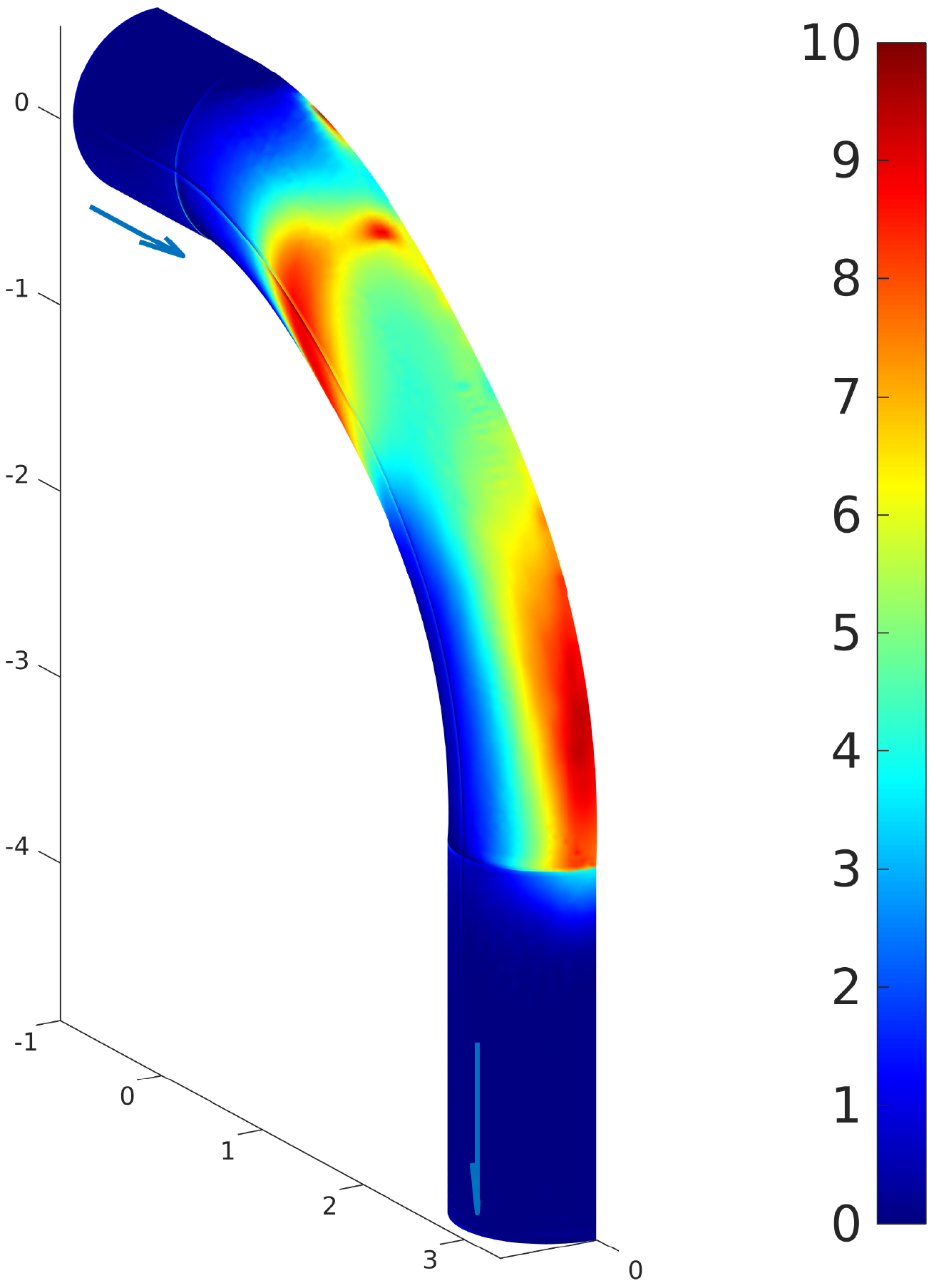}
	\subcaption{Impact angle $\dependence{\angleErosion}{\velocitySolid,\outerNormal}$ (in $^{\circ}$)}
    \label{fig:ImpactAnglesOptimized}
    \end{subfigure}
    ~
     \begin{subfigure}[t]{0.32\textwidth}
	\centering
	\includegraphics[height=5.6cm]{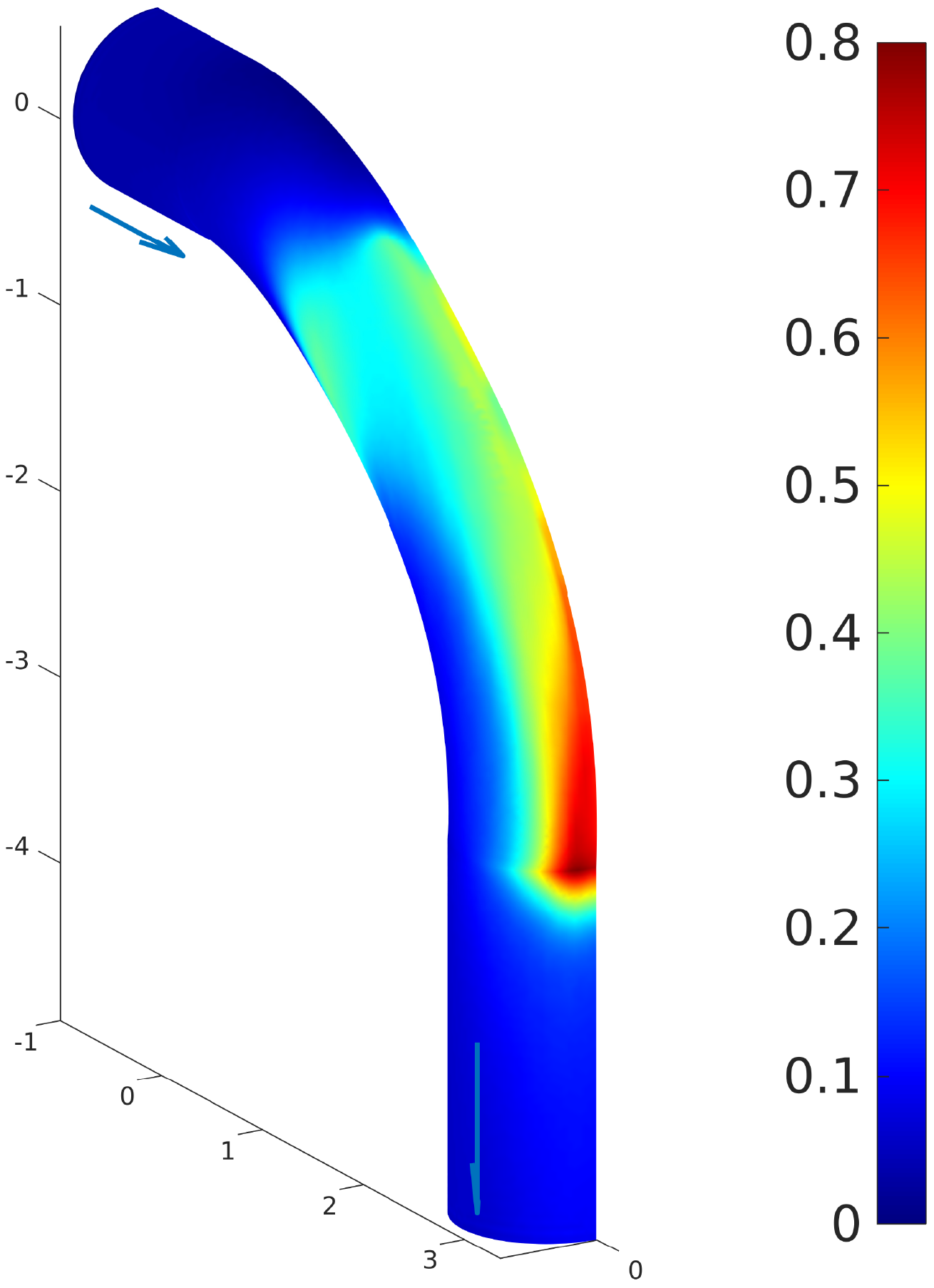}
	\subcaption{Impact velocity $\normEuclidean{\velocitySolid}$}
	\label{fig:ImpactVelocitiesOptimized}
    \end{subfigure}
    \caption{Particle impact rates, angles and velocities for the optimized geometry}
    \label{fig:ComponentsOfErosionRateOptimizedGeometry}
\end{figure}

In order to investigate this decrease of the objective more closely, 
we recall that 
$\functionalErosionLetter$ is defined in \cref{eq:ErosionRateOKA}
as the product of the local impact rate, 
the norm of the impact velocity and the impact angle dependent function $\functionAngleDependenceImpact$,
and turn our attention towards these individual contributions shown in \Cref{fig:ComponentsOfErosionRateInitialGeometry,fig:ComponentsOfErosionRateOptimizedGeometry}.
From \Cref{fig:ImpactRateInitial,fig:ImpactAnglesInitial,fig:ImpactVelocitiesInitial} it can be seen, 
that all of those three factors are relatively high in the center of the outer bended segment, which explains the relatively high erosion rates from \Cref{fig:ErosionRateInitial}. 
For the optimized geometry, the regions of the unfavorable less tangential impact angles, 
and higher impact velocities are shifted towards the lower part of the bended segment and the impact rates are distributed more evenly along $\surfaceCompDomainDeformable$, 
which can be seen in \Cref{fig:ImpactRateOptimized,fig:ImpactAnglesOptimized,fig:ImpactVelocitiesOptimized}. 
In contrast to the initial geometry however, 
the regions of high impact rates do not coincide with the regions of larger impact angles and higher impact velocities, 
which results in the reduction of the objective observed in \Cref{fig:ConvergenceHistory}.

\begin{figure}[htbp!]
  \begin{subfigure}[t]{0.49\textwidth}
	\centering
	  \begin{tikzpicture}[scale=1.00]

	  \begin{axis}[%
	  xmin=0,
	  xmax=1.4,
	  ymin=0,
	  ymax=20,
	  legend pos=south east,
	  legend cell align={left},
	  xlabel=$\stokesNumber$,
	  ylabel=
	  ]
	  \addplot [mark=*, mark options={solid}]
	    table[]{Graphics/tikz/data/ErosionRatesInitialOptimized-1.tsv};
	  \addlegendentry{$\dependence{\ErosionSquaredIntegral}{\domain_{1}}$}
	  \addplot [mark=asterisk, mark options={solid}]
	    table[]{Graphics/tikz/data/ErosionRatesInitialOptimized-2.tsv};
	  \addlegendentry{$\dependence{\ErosionSquaredIntegral}{\domain_{18}}$}
	  \end{axis}
	  \end{tikzpicture}
	\subcaption{
	    Erosion on the initial and optimized geometry
	}
	\label{fig:RelativeErosionRatesInitialAndOptimizedGeometry}
    \end{subfigure}
    ~
    \begin{subfigure}[t]{0.49\textwidth}
	\centering
	\begin{tikzpicture}[scale=1.00]
	    \begin{axis}[
		ymin=0,
		ymax=1,
		legend pos=south east,
		legend cell align={left},
		xlabel=$\stokesNumber$,
		ylabel=$\depositionFraction$
	      ]
	  
	    \addplot [mark=*, mark options={solid}]
		table[]{Graphics/tikz/data/DepositionRateOutflowToInflow_Initial.tsv}; 
		\addlegendentry{$\depositionFraction$ for $\domain_1$}
	    
	    \addplot [mark=asterisk, mark options={solid}]
		table[]{Graphics/tikz/data/DepositionRateOutflowToInflow_Optimized.tsv}; 
		\addlegendentry{$\depositionFraction$ for $\domain_{18}$}
		
		\end{axis}
	  \end{tikzpicture}
	\subcaption{
	    Impact rates $\depositionFraction$ for the initial and optimized geometry
	}
	\label{fig:DepositionRatesInitialAndOptimizedGeometry}
    \end{subfigure}  
    \caption{Erosion and impact rates on the initial and optimized geometry for different Stokes numbers}
    \label{fig:ErosionRateAnd }
\end{figure}
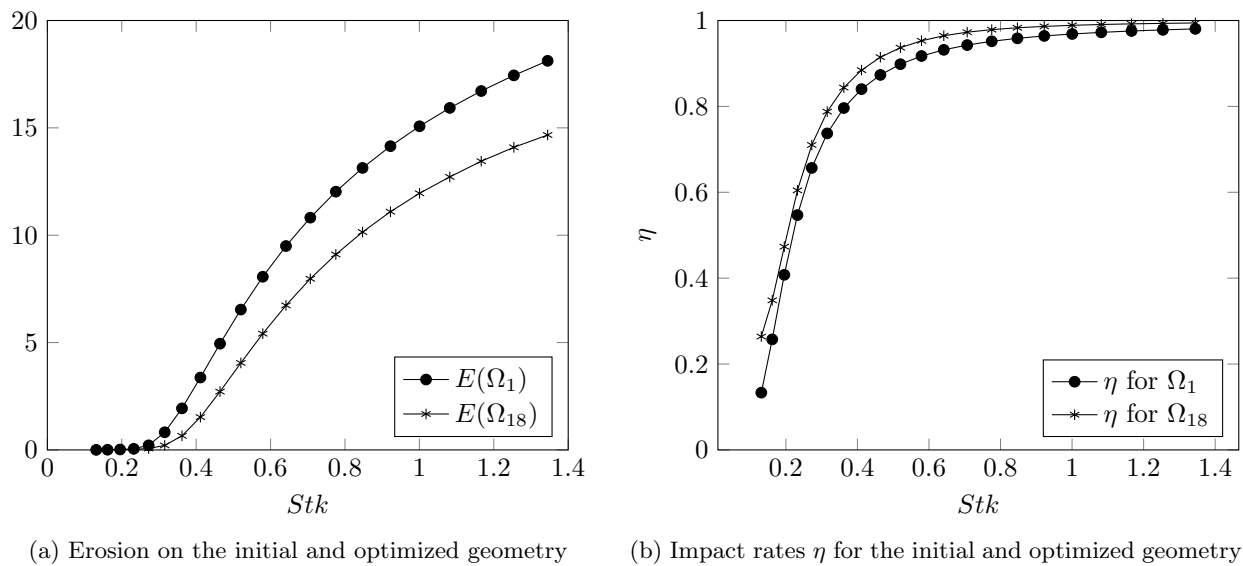

Even though we considered only a single particle species during the optimization,
we expect a continuous dependence of the erosion rate on the Stokes number in a neighborhood of the chosen value of $\stokesNumber=0.33$.
This behavior can be observed in 
\Cref{fig:RelativeErosionRatesInitialAndOptimizedGeometry},
where the integrated erosion rates
\begin{equation}
\label{eq:ErosionSquaredIntegral}
  \dependence{\ErosionSquaredIntegral}{\domain}:=
  \IntegralSurfacePartDomainPlain{\dependence{\functionalErosionLetter}{\solutionSturmApproachVolumePercentageSolid,\solutionSturmApproachVelocitySolid,\outerNormal}}{\surfaceCompDomainWall}
\end{equation}
for all of the Stokes numbers given in \Cref{tab:ParametersCollectionEfficiency}
are shown.
For all particle species with $\stokesNumber>0.2$, 
\cref{eq:ErosionSquaredIntegral}
is decreased by at least $20\%$ with respect to the corresponding value on the initial geometry.
It can be seen from \Cref{fig:DepositionRatesInitialAndOptimizedGeometry},
that these erosion rates 
are not decreased through lower impact rates, 
since they are slightly higher on the optimized geometry. 
We can therefore deduce, that 
the impact conditions for the optimized geometry
are less damaging due to smaller impact angles and reduced impact velocities
for all of the considered particle species.

%% file: Conclusion.tex
In this work we considered an optimal shape design problem, 
where the erosion due to the impact of a monodisperse particle stream
at the walls of a bended pipe segment is to be minimized.
We used a one-way coupled, volume-averaged transport model instead of Lagrangian particle tracking, 
since it allows the use of methods from PDE-constrained shape optimization,
and formally derived the adjoint equations and the shape derivative for a class of optimization problems based on this model.
We then applied the shape derivative 
within a gradient descent method for the minimization of the maximal erosion rate for  a three-dimensional test case.
We observed, that the 
erosion rate on the obtained geometry is reduced not only for the particle species, 
that is used during the optimization, 
but also for a wider range of particle diameters.
While our approach is currently restricted to laminar flow situations, 
for future work we plan to incorporate it as a coarse model within a shape optimization framework for turbulent particle erosion problems based on  space-mapping techniques \cite{Bakr2001,bandler2004space}.